\documentclass[12pt,a4paper]{article}
\usepackage{amsmath}
\usepackage{amssymb}
\usepackage{amsthm}
\usepackage{amsfonts}
\usepackage{latexsym}

\theoremstyle{plain}
\newtheorem{thm}{Theorem}[section]
\newtheorem{cor}{Corollary}[section]
\newtheorem{lem}{Lemma}[section]

\theoremstyle{definition}
\newtheorem{defn}{Definition}[section]

\theoremstyle{remark}

\newtheorem{rem}{Remark}[section]

\title{Local trace formulae and scaling asymptotics\\ in Toeplitz quantization}
\author{Roberto Paoletti\footnote{\noindent{\bf Address:}
Dipartimento di Matematica e Applicazioni, Universit\`a degli Studi
di Milano Bicocca, Via R. Cozzi 53, 20125 Milano,
Italy; {\bf e-mail}: roberto.paoletti@unimib.it }}
\date{}

\begin{document}

\maketitle

\begin{abstract}
A trace formula for Toeplitz operators was proved by Boutet de Monvel and
Guillemin in the setting of general Toeplitz structures.
Here we give a local version of this result for a class of Toeplitz
operators related to continuous groups of symmetries on quantizable compact
symplectic manifolds. The local trace formula involves certain scaling asymptotics
along the clean fixed locus of the Hamiltonian flow of the symbol, reminiscent of
the scaling asymptotics of the equivariant components of the Szeg\"{o} kernel along
the diagonal.
\end{abstract}

\section{Introduction}

A trace formula for Toeplitz operators was proved by Boutet de Monvel and Guillemin in
the setting of general Toeplitz structures \cite{bg}, following antecedents
for Laplacians \cite{c}, and more generally for positive pseudodifferential operators \cite{dg}.
The aim of the present paper is to give a local version of the trace formula in terms of suitable
scaling limits
for a special, but geometrically interesting, class of Toeplitz operators in the context of positive line bundles.

As in \cite{p-weyl}, we shall adapt the conceptual framework of \cite{z},
\cite{bsz} and \cite{sz}, where
scaling limits are studied building on the microlocal theory of the Szeg\"{o}
kernel in \cite{bs}, and pair this approach with classical arguments from
\cite{h}, \cite{dg}, \cite{gs}.
We remark that in this setting scaling limits are generally taken
with respect to the discrete
parameter indexing the isotype for the circle action on $X$; in the present situation, we
shall instead consider scaling limits with respect to the continuous auxiliary parameter in the trace
formula asymptotics.

Although quite restrictive, the class of Toeplitz operators in point is nonetheless
rather natural in geometric quantization (and algebraic geometry),
since it is related to continuous 1-parameter groups of symmetries preserving the quantization setup.
Thus the following approach
applies in particular whenever a compact Lie group acts on a quantizable compact symplectic manifold in
an Hamiltonian fashion, for in this case one can find a $G$-invariant and compatible almost complex structure,
and adopt as quantizations the spaces of sections defined by the deformation of the $\overline{\partial}$
complex introduced in \cite{bg}. The theory in \cite{sz} provides a generalization to this context
of the microlocal description of the Szeg\"{o} kernel in \cite{bs}.
For ease of exposition, we shall confine ourselves to the more familiar complex projective case.

Thus let $M$ be a d-dimensional complex projective manifold, and let $A$ be an ample
line bundle on it. Let $h$ be an Hermitian metric on $A$, and $\nabla$ be the
unique covariant derivative compatible with the holomorphic and Hermitian
structures. Assume, as we may, that $\nabla$ has curvature $\Theta=-2i\omega$, where
$\omega$ is a K\"{a}hler form. Endowed with the volume form
$dV_M=:(1/\mathrm{d}!)\,\omega^{\wedge \mathrm{d}}$, $M$ has total volume
$\mathrm{vol}(M)=\left(\pi^{\mathrm{d}}/\mathrm{d}!\right)\,\int_Mc_1(A)^\mathrm{d}$.

Our focus here is on Hamiltonians generating 1-parameter flows of
holomorphic symplectomorphisms of $M$. More precisely, given a real function $f\in \mathcal{C}^\infty(M)$ let $\upsilon_f$ be its
Hamiltonian vector field in the symplectic structure $2\omega$; since $M$ is compact,  $\upsilon_f$ generates
a 1-parameter group of Hamiltonian symplectomorphisms
$
\phi^M:\tau\in \mathbb{R}\mapsto \phi^M_\tau\in \mathrm{Symp}(M)$.

\begin{defn}
We shall call $f$ a  \textit{compatible Hamiltonian} if every $\phi_\tau^M:M\rightarrow M$ is holomorphic.
\end{defn}

For instance,
if $M=\mathbb{CP}^1$ and $A$ is the hyperplane bundle endowed with the Fubini-Study metric,
then $f\big([z_0:z_1]\big)=:\left(k|z_0|^2+l|z_1|^2\right)/\left(|z_0|^2+|z_1|^2\right)$
is compatible for any pair of integers $k,l$.
More generally, let $G$ be a compact Lie group with Lie algebra
$\mathfrak{g}$, and let $G\times M\rightarrow M$ be an holomorphic Hamiltonian action.
Any $\xi\in \mathfrak{g}$ induces a vector field
$\xi_M$ on $M$.
If $\Phi:M\rightarrow \mathfrak{g}^*$ is the moment map and  $f=:\langle \Phi,\xi\rangle$,
then $\upsilon_f=\xi_M$ and $f$ is compatible.

If $f$ is compatible, $\phi_\tau^M$ induces a 1-parameter
group of
unitary automorphisms
$
\psi_\tau^{(k)}:H^0\left(M,A^{\otimes k}\right)\rightarrow H^0\left(M,A^{\otimes k}\right)$ for
$k=0,1,2,\ldots$;
here $H^0\left(M,A^{\otimes k}\right)$ is the space of global holomorphic sections of
$A^{\otimes k}$, endowed with the natural Hermitian product induced by $h$ and $dV_M$.
In fact, $\phi^M$ lifts to a 1-parameter group of holomorphic bundle automorphism $\phi^A_\tau:A\rightarrow A$, and
$\psi^{(1)}_\tau(s)=\phi^A_\tau\circ s\circ \phi^M_{-\tau}$ for $s\in H^0(M,A)$; similarly for $k\ge 2$.

Let $\mathfrak{H}(A)=:\bigoplus_{k}H^0\left(M,A^{\otimes k}\right)$
be the Hilbert space direct sum;
then $
\psi_\tau=:\oplus _k\psi_\tau^{(k)}:\mathfrak{H}(A)\rightarrow
\mathfrak{H}(A)$ is a unitary isomorphism.
If $f>0$, the trace of $\psi_\tau$ is a tempered distribution on $\mathbb{R}$, and its singular support is a
set of periods for an appropriate Hamiltonian flow; the trace formula describes
its singularities at each isolated period (the flow in point is not the flow of $\xi_f$ on $M$,
but is closely related to it). These singularities are encapsulated in certain asymptotic expansions,
whose coefficients relate to the dynamics of the closed trajectories of the given period; in particular, the leading term
is given by Poincar\'{e}
map data.

To make this more precise, it is convenient to lift the problem to the unit circle bundle
$A^*\supseteq X\stackrel{\pi}{\rightarrow}M$.
The connection form $\alpha$ on $X$ is then a contact form, and $d\mu_X=:(1/2\pi)\,\alpha\wedge \pi^*(dV_X)$
is a volume form. We shall henceforth identify functions, densities and half-densities on $X$.

The Hamiltonian vector field $\upsilon_f$ lifts to a contact vector field
$\widetilde{\upsilon}_f$ on $X$; let  $\phi_\tau^X:X\rightarrow X$ be the
associated 1-parameter group of contactomorphisms.
Pull-back by $\phi_\tau^X$,  $\left(\phi_{\tau}^X\right)^*:L^2(X)\rightarrow L^2(X)$, is a unitary isomorphism.

In addition, the hypothesis that $\phi_\tau^M$ be holomorphic
implies that $\left(\phi_{\tau}^X\right)^*$
leaves the Hardy space $H(X)\subseteq L^2(X)$ invariant.
Under the natural unitary isomorphism $H(X)\cong\mathfrak{H}(A)$, $\psi_\tau$ corresponds to
$\left(\phi_{-\tau}^X\right)^*$.

Let $U_H(\tau):H(X)\rightarrow H(X)$
be the unitary operator induced by $\left(\phi_{-\tau}^X\right)^*$, and let us extend $U_H(\tau)$
to $L^2(X)$ by declaring it to vanish on the orthocomplement $H(X)^\perp$. In other words,
$U_H(\tau)=\left(\phi_{-\tau}^X\right)^*\circ \Pi$, where $\Pi:L^2(X)\rightarrow H(X)$
is the orthogonal projector. Basic results on wave fronts imply
that $U_H(\tau)$ extends to a continuous operator $\mathcal{D}'(X)\rightarrow \mathcal{D}'(X)$ \cite{d};
we shall also denote by
$U_H(\tau)\in \mathcal{D}'(X\times X)$ its distributional kernel.
As a differential operator on $X$, $i\upsilon_f$ leaves
$H(X)\cap \mathcal{C}^\infty(X)$ invariant; therefore, its restriction is the first order
self-adjoint Toeplitz operator $T_f=:i\upsilon_f\circ \Pi$, which has a positive symbol.
Let $\lambda_1\le \lambda_2\le \ldots$ be the eigenvalues of $T_f$ acting on $H(X)$ \cite{bg}, so that
$e^{i\lambda_j\tau}$ are the eigenvalues of $U_H(\tau)$.
The trace in point is the tempered distribution $\sum_je^{i\lambda_j\tau}$ on $\mathbb{R}$.

The asymptotic estimates describing the singularities of $\sum_je^{i\lambda_j\tau}$
involve certain functions $S_\chi\in \mathcal{C}^\infty(X\times X)$.
More precisely, for every $\chi\in \mathcal{C}^\infty_0(\mathbb{R})$ the averaged operator
$$
S_\chi=:\int_{-\infty}^{+\infty}\chi(\tau)\,U_H(\tau)\,d\tau
$$
is smoothing, and thus (with abuse of language) has kernel
$S_\chi\in \mathcal{C}^\infty(X\times X)$. In particular,
$$
\left<\sum_je^{i\lambda_j\tau},\chi\right>=\mathrm{trace}(S_\chi)=\int_XS_\chi(x,x)\,dV_X(x).
$$

\begin{defn}
$\mathrm{Per}_X(f)\subseteq \mathbb{R}$ is the set of periods of $\phi^X$. Thus
$\tau\in \mathrm{Per}_X(f)$ if and only if
$\phi_\tau^X(x)=x$ for some $x\in X$.
\end{defn}

If $\tau\in \mathrm{Per}_X(f)$, then $\tau$ is also a period of $\phi^M$, but the converse needn't be true.
However, let $\mathrm{Fix}\left(\phi_\tau^M\right)\subseteq M$ and
$\mathrm{Fix}\left(\phi_\tau^X\right)\subseteq X$ be the fixed loci of $\phi_\tau^M$ and $\phi_\tau^X$, respectively;
if $\tau\in \mathrm{Per}_X(f)$, then by $S^1$-invariance $\phi_\tau^X$ is an $S^1$-bundle over a union of connected
components of $\mathrm{Fix}\left(\phi_\tau^M\right)$.

Now $\mathrm{Per}_X(f)$ contains the singular support of $\sum_je^{i\lambda_j\tau}$.
Suppose $\tau_0$ is an isolated point of $\sum_je^{i\lambda_j\tau}$, and $\chi$ is a bump function
with $\chi(\tau_0)=1$ and supported in a small open neighborhood
of $\tau_0$. Then $\chi(\tau)\,\sum_je^{i\lambda_j\tau}$ equals $\sum_je^{i\lambda_j\tau}$ near
$\tau_0$, and is compactly supported and non-singular away from $\tau_0$;
therefore,
the singularity of $\sum_je^{i\lambda_j\tau}$ at $\tau_0$
is characterized by the asymptotics of the Fourier transform of
$\chi(\tau)\,\sum_je^{i\lambda_j\tau}$, \textit{viz.} of the trace of
$S_{\chi\,e^{-i\lambda (\cdot)}}$ as $\lambda\rightarrow \infty$.

Here we shall localize the problem, and study the pointwise asymptotics of the kernel of
$S_{\chi\,e^{-i\lambda (\cdot)}}$; this amounts to considering certain scaling limits in
the neighborhood of the fixed locus of $\phi^X_{\tau_0}$.
A global trace formula follows by integration.

The following estimates are phrased in terms of local Heisenberg coordinates on $X$
\cite{sz}. Having fixed a system of Heisenberg local coordinates
$(\theta,\mathbf{v})$ centered at $x\in X$, we shall follow \cite{sz} and write
$x+(\theta,\mathbf{v})$ for the point with local coordinates $(\theta,\mathbf{v})$, and
$x+\mathbf{v}$ for $x+(0,\mathbf{v})$; we have $x+(\theta,\mathbf{v})=
\mathfrak{r}_{e^{i\theta}}(x+\mathbf{v})$, where $\mathfrak{r}:S^1\times X\rightarrow X$
is the circle action.
A choice of Heisenberg local coordinates on $X$ centered at $x_0$ includes a choice of
\textit{preferred} (not necessarily holomorphic)
local coordinates on $M$ centered at $m_0=:\pi(x_0)$, whence
of a unitary isomorphism $T_{m_0}M\cong \mathbb{C}^\mathrm{d}$. With this implicit,
the expression $x+\mathbf{v}$
may be used for either $\mathbf{v}\in T_{m_0}M$ or $\mathbf{v}\in \mathbb{C}^\mathrm{d}$.

Before stating the Theorem, let us introduce two further pieces of notation.

First, following  \cite{sz}, we shall denote by $\psi_2:\mathbb{C}^\mathrm{d}\times \mathbb{C}^\mathrm{d}\rightarrow \mathbb{C}$ the
smooth function:
$$
\psi_2(\mathbf{u},\mathbf{w})=:i\,\Im\left(\mathbf{u}\cdot _c\overline{\mathbf{w}}\right)
-\frac 12\,\|\mathbf{u}-\mathbf{v}\|^2,
$$
where $\mathbf{a}\cdot _c\mathbf{b}\in \mathbb{C}$ is the standard scalar product of vectors $\mathbf{a},\mathbf{b}\in
\mathbb{C}^\mathrm{d}$, so that $\mathbf{a}\cdot _c\overline{\mathbf{b}}$ is their standard Hermitian product
($\mathbf{a}\cdot\mathbf{b}=:\Re(\mathbf{a}\cdot _c\mathbf{b})\in \mathbb{R}$ will denote their standard Euclidean product
as vectors in $\mathbb{R}^{2\mathrm{d}}$).

If $f$ is compatible, and $\tau_0\in \mathrm{Per}_X(f)$ has a clean fixed locus,
then $\mathrm{Fix}\left(\phi_{\tau_0}^M\right)\subseteq M$ is a complex submanifold, hence
for any $x_0\in \mathrm{Fix}\left(\phi_{\tau_0}^X\right)$ its tangent space at $\pi(x_0)$, $T_{\pi(x_0)}\left(\mathrm{Fix}\left(\phi_{\tau_0}^M\right)\right)$, is a complex vector subspace of the
tangent space $T_{\pi(x_0)}M$ to $M$.
Let us denote by $T_{\pi(x_0)}\left(\mathrm{Fix}\left(\phi_{\tau_0}^M\right)\right)^\perp\subseteq T_{\pi(x_0)}M$
the orthocomplement.

\begin{thm}
\label{thm:main}
Assume that $f\in \mathcal{C}^\infty(M)$ is positive and compatible.
Let $\widetilde{\upsilon}_f$ be the contact vector field on $X$ induced by
$f$, and let $\lambda_1\le \lambda_2\le \cdots$ be the eigenvalues of the
Hermitian operator $i\widetilde{\upsilon}_f$ acting on the Hardy space $H(X)$.
Then the singular support of the tempered distribution $\sum_je^{i\lambda_j\tau}$ is contained in $\mathrm{Per}_X(f)$.
Furthermore, assume that $\tau_0\in \mathrm{Per}_X(f)$ is an isolated period with clean
fixed locus $\mathrm{Fix}\left(\phi_{\tau_0}^X\right)\subseteq X$. Then there exists $\epsilon>0$ such that for
all $\chi\in \mathcal{C}^\infty\big((\tau_0-\epsilon,\tau_0+\epsilon)\big)$ and $C>0$ the following holds.
\begin{enumerate}
\item $
  S_{\chi\,e^{-i\lambda\,(\cdot) }}(x,x)=O\left(\lambda^{-\infty}\right)
  $ uniformly on $X\times X$ as $\lambda\rightarrow -\infty$;
  \item $
  S_{\chi\,e^{-i\lambda\,(\cdot) }}(x,x)=O\left(\lambda^{-\infty}\right)
  $
  uniformly for $\mathrm{dist}_X\Big(x,\mathrm{Fix}\left(\phi_{\tau_0}^X\right)\Big)\ge C\,
  \lambda^{-7/18}$ as $\lambda\rightarrow +\infty$.
    \item Uniformly in
$x_0\in \mathrm{Fix}\left(\phi_{\tau_0}^X\right)$ and
$
\mathbf{u}\in
T_{\pi(x_0)}M
$
with
\begin{center}
$\|\mathbf{u}\|\le C\lambda^{1/9}$ and
$\mathbf{u}\in T_{\pi(x_0)}\left(\mathrm{Fix}\left(\phi_{\tau_0}^M\right)\right)^\perp$,
\end{center}
as $\lambda\rightarrow +\infty$ we have an asymptotic expansion
\begin{eqnarray}
\label{eqn:thm-main-expansion}
\lefteqn{S_{\chi\,e^{-i\lambda(\cdot)}}\left(x_0+\frac{\mathbf{u}}{\sqrt{\lambda}},x_0+\frac{\mathbf{u}}{\sqrt{\lambda}}\right)}\\
&\sim&
\frac{2\pi\,e^{-i\lambda\tau_0}}{f(m_0)^{\mathrm{d}+1}}\left(\frac{\lambda}{\pi}\right)^{\mathrm{d}}\,
e^{f(m_0)^{-1}\cdot\psi_2\big(d_{m_0}\phi^M_{-\tau_0}(\mathbf{u}),\mathbf{u}\big)}\,\chi(\tau_0)\cdot \left[1+\sum _{j=1}^{+\infty}\lambda^{-j/2}G_j(x_0,\mathbf{u})\right].\nonumber
\end{eqnarray}
where $G_j(x_0,\mathbf{u})$'s are polynomials in $\mathbf{u}$ depending smoothly
on $x_0$. More precisely, in the given range the $N$-th remainder is uniformly
$O\left(\lambda^{-aN}\right)$ for some $a>0$.
\item Let us write
$$
S_{\chi\,e^{-i\lambda(\cdot)}}\left(x_0+\frac{\mathbf{u}}{\sqrt{\lambda}},x_0+\frac{\mathbf{u}}{\sqrt{\lambda}}\right)
=\mathfrak{E}_\lambda(\mathbf{u})+\mathfrak{O}_\lambda(\mathbf{u}),
$$
where $\mathfrak{E}_\lambda$ and $\mathfrak{O}_\lambda$ are even and odd functions of $\mathbf{u}$, respectively.
Then the asymptotic expansion of $\mathfrak{E}_\lambda(\mathbf{u})$
as $\lambda\rightarrow +\infty$ is obtained from (\ref{eqn:thm-main-expansion})
by collecting all integer powers of $\lambda$; similarly, the asymptotic expansion of $\mathfrak{O}_\lambda(\mathbf{u})$
is obtained by collecting all fractional (non-integer) powers of
$\lambda$.
\end{enumerate}

\end{thm}

We obtain a global trace formula by integration. To state this, we need a further piece of
notation. In the hypothesis of the Theorem, $X_{\tau_0}=:\mathrm{Fix}\left(\phi_{\tau_0}^X\right)$
is a union of connected components $X_{\tau_0,j}$, $1\le j\le N_{\tau_0}$, of real dimension
$\dim\left(X_{\tau_0,j}\right)=2\mathrm{f}_j+1$; thus, $\mathrm{f}_j$ is the complex dimension of the connected
component $M_{\tau_0,j}=:\pi\left(X_{0j}\right)$ of $M_{\tau_0}=:\mathrm{Fix}\left(\phi_{\tau_0}^M\right)$.

If
$m_0\in M_{\tau_0, j}$, let $N_{m_0}=:\left(T_{m_0}M_{\tau_0, j}\right)^\perp\subseteq T_{m_0}M$ be the normal space to $M_{\tau_0, j}$ at $m_0$.
As the fixed locus is clean, the differential $d_{m_0}\phi^M_{\tau_0}$
restricts to a unitary isomorphism
$\left.d_{m_0}\phi^M_{\tau_0}\right|_{N_{m_0}}:N_{m_0}\rightarrow N_{m_0}$ such that
$\mathrm{id}_{N_{m_0}}-\left.d_{m_0}\phi^M_{\tau_0}\right|_{N_{m_0}}$ is an isomorphism.
The determinant of the latter isomorphism only depends on $j$. We may set
$$
c(\tau_0,j)=:\det\left(\mathrm{id}_{N_{m_0}}-\left.d_{m_0}\phi^M_{-\tau_0}\right|_{N_{m_0}}\right)\,\,\,\,\,\,\,\,
\left(m_0\in M_{\tau_0, j}\right).
$$

Let $dV_{M_{\tau_0 j}}$ be the volume form on $M_{\tau_0, j}$ given by restriction
of $\omega^{\wedge\mathrm{f}_0}/\mathrm{f}_0!$.

\begin{cor}
\label{cor:global-trace-formula}
In the hypothesis of Theorem \ref{thm:main}, and with the above notation, the following holds.
As $\lambda\rightarrow -\infty$ we have
$$
\int_XS_{\chi\,e^{-i\lambda(\cdot)}}(x,x)\,dV_X(x)=O\left(\lambda^{-\infty}\right);
$$
as $\lambda\rightarrow +\infty$, on the other hand, we have
$$
\int_XS_{\chi\,e^{-i\lambda(\cdot)}}(x,x)\,dV_X(x)=\sum _{j=1}^{N_{\tau_0}}I_j(\chi,\lambda),
$$
where each $I_j$ admits an asymptotic expansion
\begin{eqnarray*}
I_j(\chi,\lambda)
&\sim&2\pi\,e^{-i\lambda\tau_0}\,\left(\frac{\lambda}{\pi}\right)^{\mathrm{f}_j}\,
\frac{\chi(\tau_0)}{c(\tau_0,j)}\cdot\left(\int_{M_{\tau_0 j}}\frac{1}{f(m)^{\mathrm{f}_j+1}}\,
dV_{M_{\tau_0 j}}(m)\right)\\
&&\cdot \left(1+\sum _{j=1}^{+\infty}\lambda^{-j}p_j\right).
\end{eqnarray*}
\end{cor}

The following remarks are in order.

The Corollary expresses
$
\Gamma(\lambda)=:\mathrm{trace}\left(S_{\chi e^{-i\lambda(\cdot)}}\right)
$
as an asymptotic expansion with general term a multiple of $e^{-i\lambda\tau_0}\lambda_+^a$,
where $a$ is a descending sequence of integers. Since $\Gamma(\lambda)$ is the Fourier transform of
$\chi(\tau)\,\mathrm{trace}\big(U(\tau)\big)$, the latter is given by an asymptotic expansion in decreasing
powers of $(\tau-\tau_0+i\,0)^{-1}$.
More precisely, there is one such expansion for each $X_{0j}$, and the
leading term of the
$j$-th expansion is a scalar multiple of $\left(\tau-\tau_0+i \,0^+\right)^{-(\mathrm{f}_j+1)}$ (cfr \cite{bg}).

The local asymptotics related to the singularity at $\tau=0$ have already been studied
in a more general setting in \cite{p-weyl}.

With obvious changes, the proof of Theorem \ref{thm:main} applies to
more general scaled pairs of the form $(x_0+\mathbf{w}/\sqrt{\lambda},x_0+\mathbf{u}/\sqrt{\lambda})$.

For a slight generalization of Theorem \ref{thm:main}, one may consider any
$S^1$-invariant first order self-adjoint Toeplitz operator $T=\Pi\circ i\widetilde{\upsilon}_f+T'$,
with $T'$ of degree $0$; the proof carries over to this situation,
but the leading order term in the asymptotic
expansion will now depend on the subprincipal symbol of $T$ through an oscillatory factor.

\section{Preliminaries.}

The cotangent bundle $T^*X$ is endowed with the standard symplectic structure $\omega_{\mathrm{stan}}$. Explicitly, suppose given
local coordinates $\mathbf{q}$ on an open subset $U\subseteq X$, and let $(\mathbf{q},\mathbf{p})$ the corresponding
 local coordinates on $TU\subseteq TX$, so that $(\mathbf{q}_0,\mathbf{p}_0)$ corresponds to the tangent vector
 $\mathbf{p}_0\cdot\left.(\partial/\partial \mathbf{q})\right|_{\mathbf{q}_0}$; then $\omega_{\mathrm{stan}}=
 d\mathbf{p}\wedge d\mathbf{q}$ on $TU$.

Let $\Sigma=:\big\{(x,r\alpha_x):x\in X,\,r>0\big\}$ be the closed symplectic cone in
$T^*X\setminus\{0\}$ generated by the connection 1-form.
Thus $\Sigma\cong X\times \mathbb{R}_+$ canonically.
Let $\theta$ be the \lq circle\rq\, coordinate on $X$ (locally defined) and let $r$ be
the \lq cone\rq\, coordinate
on $\Sigma$; then
the restriction of $\omega_{\mathrm{stan}}$ to $\Sigma$ is
$\omega_\Sigma=2r\omega+dr\wedge d\theta$
(symbols of pull-back are omitted).

We adopt the convention that the Hamiltonian vector field
$\upsilon_g$ of a sooth function $g$ with respect to a symplectic structure
$\Omega$ is defined by $\Omega\big(\iota(\upsilon_g),u)=dg(u)$ for any
tangent vector $u$.
The Hamiltonian vector field on $\Sigma$ of
$\widetilde{f}=:rf\in \mathcal{C}^\infty(\Sigma)$ is then $\widetilde{\upsilon}_f$,
where $\widetilde{\upsilon}_f=\upsilon_f^\sharp-f\,(\partial/\partial \theta)$ is the contact lift of
$\upsilon_f$.
Here $\upsilon_f^\sharp$ is the horizontal lift of $\upsilon_f$, and
$\partial/\partial \theta$ is the infinitesimal generator of the $S^1$-action
(with some ambiguity,
$\widetilde{\upsilon}_f$ denotes both the contact lift of $\upsilon_f$ to $X$ and its further Hamiltonian lift
to $\Sigma$; strictly speaking, the latter is $\big( \widetilde{\upsilon}_f,0\big)$).

The 1-parameter group of symplectomorphisms
$\phi^\Sigma_\tau:\Sigma\rightarrow \Sigma$ generated by $\widetilde{\upsilon}_f$ is
$
\phi^\Sigma_\tau\big(x,r\alpha_x\big)=:\left(\phi^X_{\tau}(x),r\alpha_{\phi^X_{\tau}(x)}\right)$.

Viewed as a differential operator on $X$, by our hypothesis $\widetilde{\upsilon}_f$ commutes with $\Pi$ and is elliptic
on $\Sigma$, since its symbol there is $-i\widetilde{f}$; therefore, it remains elliptic in a conic neighborhood
of $\Sigma$ in $T^*X\setminus\{0\}$.
Hence $i\widetilde{\upsilon}_f$ restricted to $H(X)$ is an elliptic self-adjoint Toeplitz operator $T_f$ of the first order,
with symbol
$\sigma_{T_f}=\widetilde{f}:\Sigma\rightarrow \mathbb{R}$.
By the proof of Lemma 12.2 of \cite{bg}, there exists a self-adjoint first
order pseudodifferential operator $Q$ on $X$ with everywhere positive principal symbol $q>0$, which
commutes with $\Pi$, equals $i\,\widetilde{\upsilon}_f$ on $H(X)$ (thus $T_f=Q\circ \Pi$),
and is such that $Q-i\,\widetilde{\upsilon}_f$ is microlocally
smoothing on a conic neighbourhood of $\Sigma$.

The operator $U_H(\tau)$, on the other hand, may be written
$
U_H(\tau)=U(\tau)\circ \Pi
$,
where $U(\tau)=e^{i\tau Q}$; by the previous discussion, its wave front is
$$
\mathrm{WF}\big(U_H(\tau)\big)=
\Big\{(x,r\alpha_x,y,-r\alpha_y)\,:\,x\in X,y=\phi_{-\tau}^X(x),\,r>0\Big\},
$$
and so the singular support is
$
\mathrm{sing.\,supp}\big(U_H(\tau)\big)=
\mathrm{graph}\left(\phi_{-\tau}^X\right)$.

Proposition 12.4 of \cite{bg} shows by a functorial description of $\sum_je^{i\lambda_j\tau}$
that its singular support is contained in the set of periods of $\phi_\tau^\Sigma$, which
is the same as the set of periods of $\phi_\tau^X$.

Now suppose $x\not\in \mathrm{Fix}\left(\phi_{\tau_0}^X\right)$, that is,
$(x,x)\not\in \mathrm{graph}\left(\phi_{\tau_0}^X\right)$. By continuity, for some
$\epsilon>0$ we have
$(x,x)\not\in \mathrm{graph}\left(\phi_{-\tau}^X\right)$ for any $\tau\in (\tau_0-2\epsilon,\tau_0+2\epsilon)$.
As a distribution on $\mathbb{R}\times X\times X$, $U_H$ is $\mathcal{C}^\infty$
on $(\tau_0-\epsilon,\tau_0+\epsilon)\times X'\times X'$, where $X'\subseteq X$ is an appropriate open neighborhood of
$X$. If $\chi\in \mathcal{C}^\infty_0\big((\tau_0-\epsilon,\tau_0+\epsilon)\big)$, then
$\chi(\tau)\,U_H\left(\tau,x',x''\right)$ is $\mathcal{C}^\infty$ on $\mathbb{R}\times X'\times X'$; hence its
Fourier transform in $\tau$ is of rapid decrease in $\lambda$ uniformly in
$\left(x',x''\right)\in X''\times X''$, where $X''\Subset X'$ is an open neighborhood of $x$. In particular,
$$
S_{\chi\,e^{-i\lambda(\cdot)}}(x,x)=O\left(\lambda^{-\infty}\right)
$$
as $\lambda\rightarrow \infty$, uniformly in $x$ with $\mathrm{dist}_X\left(x,\mathrm{Fix}\left(\phi_{\tau_0}^X\right)\right)\ge \delta$
for any fixed $\delta>0$ (here $\mathrm{dist}_X$ is the Riemannian distance function on $X$), for any $\chi$ with sufficiently small
compact support near $\tau_0$.

\section{Proof of Theorem \ref{thm:main}.}

We shall adapt the approach used in \cite{p-weyl}
for the singularity at $\tau=0$ to the case of a general period (and thus occasionally refer to arguments in \cite{p-weyl}).

When working locally on $X$, looking for asymptotic expansions we may replace $U(\tau)$ and $\Pi$ by their respective representations as Fourier integral operators.
Thus we shall write \cite{dg}
$$
U(\tau)(x,y)=\frac{1}{(2\pi)^{2\mathrm{d}+1}}\,\int_{\mathbb{R}^{2\mathrm{d}+1}}e^{i\big[\varphi(\tau,x,\eta)-y\cdot \eta\big]}\,
a(\tau,x,y,\eta)\,d\eta,
$$
where $a(\tau,\cdot,\cdot)\in S^0_{\mathrm{cl}}$, and
the phase $\varphi(\tau,\cdot,\cdot)$ is a local generating function for the 1-parameter group
$\phi^{T^*X}_{-\tau}$ of homogeneous Hamiltonian symplectomorphism of $T^*X$
generated by $-q$. Working with $\tau$ close to $\tau_0$, we change variable
$\tau\rightsquigarrow \tau+\tau_0$, so that
\begin{equation}
\label{eqn:hamilton-jacobi}
\varphi\big(\tau+\tau_0,x,\eta\big)=\varphi(\tau_0,x,\eta)+\tau\,q\big(x,d_x\varphi(\tau_0,x,\eta)\big)+O\left(\tau^2\right)\,\|\eta\|.
\end{equation}
On the other hand, $\Pi$ is an FIO with complex phase of the form
$$
\Pi(x,y)=\int_0^{+\infty}e^{it\psi(x,y)}\,s(t,x,y)\,dt,
$$
where $s$ is a classical symbol $s(t,x,y)\sim \sum_{j\ge 0}s_j(x,y)\,t^{\mathrm{d}-j}$, and the Taylor expansion of
$\psi$ along the diagonal is determined by the K\"{a}hler metric \cite{bs}.

Let $\sim$ stand for \lq equal asymptotics as $\lambda\rightarrow \infty$\rq.
Given any $x\in X$, let $\varrho $ be a positive smooth cut-off function on $X$, identically equal to $1$ near $x$,
and with support contained in a sufficiently open neighborhood $X_1\Subset X$;
arguing as in \cite{p-weyl},
 Lemma 2.1 we
get
\begin{eqnarray}
\label{eqn:1st-oscillatory-expression}
\lefteqn{S_{\chi\,e^{-i\lambda(\cdot)}}(x,x)}\\
&\sim&\frac{1}{(2\pi)^{2\mathrm{d}+1}}\,\int_{X_1}\int _{-\epsilon}^{\epsilon}e^{-i\lambda(\tau+\tau_0)}
\chi(\tau+\tau_0)\,\varrho(z)U(\tau+\tau_0)(x,z)\,\Pi(z,x)\,d\mu_X(z)\,d\tau \nonumber\\
&\sim&\frac{1}{(2\pi)^{2\mathrm{d}+1}}\,\int_{X_1}\left[\int_{\mathbb{R}^{2\mathrm{d}+1}}\int _{-\epsilon}^{\epsilon}e^{i\Psi(\tau,x,z,\eta,\lambda)}
B(\tau,x,z,\eta)\,d\tau\,d\eta\right]\cdot\Pi(z,x)\,d\mu_X(z),\nonumber
\end{eqnarray}
where
$$
B(\tau,x,z,\eta)=:\chi(\tau+\tau_0)\,\varrho(z)\,a(\tau+\tau_0,x,z,\eta),
$$
and in view of (\ref{eqn:hamilton-jacobi})
\begin{eqnarray*}
\lefteqn{\Psi(\tau,x,z,\eta,\lambda)=:\varphi(\tau+\tau_0,x,\eta)-z\cdot \eta
-\lambda(\tau+\tau_0)}\\
&=&\varphi(\tau_0,x,\eta)+\tau\,q\big(x,d_x\varphi(\tau_0,x,\eta)\big)+O\left(\tau^2\right)\,\|\eta\|
-y\cdot \eta
-\lambda(\tau+\tau_0).
\end{eqnarray*}
Now as $\lambda\rightarrow -\infty$ since $q$ is elliptic and positive for sufficiently small $\epsilon$
we have
$$
\partial _\tau\Psi=q\big(x,d_x\varphi(\tau_0,x,\eta)\big)+O\left(\tau\right)\,\|\eta\|-\lambda\ge C_1\,\big(\|\eta\|+|\lambda|\big)
$$
for some fixed $C_1>0$;
by repeated integration by parts in $d\tau$ we conclude
$S_{\chi\,e^{-i\lambda(\cdot)}}(x,x)=O\left(\lambda^{-\infty}\right)$ uniformly in $x\in X$ as $\lambda\rightarrow
-\infty$.

Next we consider the asymptotics for $\lambda\rightarrow +\infty$. By the discussion in \S 1, we are reduced to considering
the problem in the neighborhood of $\mathrm{Fix}\left(\phi_{\tau_0}^X\right)$, so we fix $x_0\in \mathrm{Fix}\left(\phi_{\tau_0}^X\right)$ and a suitably small open neighborhood $X_1\subseteq X$ of $x_0$, and consider the
asymptotics for $x\in X_1$. Let there be given on $X_1$ a system of Heisenberg local coordinates
$(\theta,\mathbf{ v}):X_1\rightarrow (-\delta,\delta)\times B_{2\mathrm{d}}(\mathbf{0},\delta)$ centered at $x_0$
for some $\delta>0$, where $B_{2\mathrm{d}}(\mathbf{0},\delta)\subseteq \mathbb{R}^{2\mathrm{d}}\cong
\mathbb{C}^\mathrm{d}$ is the open ball of radius $\delta>0$ centered at the origin. In particular,
in the associated cotangent coordinates $(x_0,\alpha_{x_0})\in T^*X$ corresponds to $\big((0,\mathbf{0}),(1,\mathbf{0})\big)$.

To begin with, we only lose a smoothing term, hence a negligible contribution to the asymptotics, if
we multiply the amplitude $a$ of $U$ by a radial function $b(\eta)$ identically equal to $0$ near the
origin and to $1$ for $\|\eta\|\gg 0$. We may then assume without loss that $a$ vanishes identically near
the origin as a function of $\eta$.

Let $S_1,S_2\subseteq S^{2\mathrm{d}}$ be open subsets
covering $S^{2\mathrm{d}}$ with $(1,\mathbf{0})\not\in \overline{S_2}$, and
let $\gamma_1+\gamma_2=1$ be a partition of unity on $ S^{2\mathrm{d}}$ subordinate to the open cover
$\{S_1,S_2\}$. Also, for $j=1,2$ let $U^{(j)}(\tau)$ be defined as $U(\tau)$, but with the amplitude $a$
multiplied by $\varrho(z)\,\gamma_j\big(\eta/\|\eta\|\big)$, and let $S_{\chi\,e^{-i\lambda(\cdot)}}^{(j)}$
be defined as $S_{\chi\,e^{-i\lambda(\cdot)}}$ with $U^{(j)}(\tau)$ in place of $U(\tau)$.
Then clearly
\begin{eqnarray}
\label{eqn:0th-oscillatory-expression-romega}
S_{\chi\,e^{-i\lambda(\cdot)}}(x,x)\sim\sum_{j=1}^2S_{\chi\,e^{-i\lambda(\cdot)}}^{(j)}(x,x).
\end{eqnarray}

Now we claim that, perhaps after replacing $X_1$ with a smaller open neighborhood and $\epsilon$ with a smaller
positive real number, the second summand on the right hand side of $(\ref{eqn:0th-oscillatory-expression-romega})$
gives a negligible contribution to the asymptotics, uniformly in $x\in X_1$.

To see this, consider the operator giving rise to the second summand:
\begin{eqnarray}
\label{eqn:2nd-summand}
S_{\chi\,e^{-i\lambda(\cdot)}}^{(j)}=
\int_{-\epsilon}^{+\epsilon}e^{-i\lambda \tau}\,\Big(U^{(2)}(\tau+\tau_0)\circ \Pi\Big)\,d\tau.
\end{eqnarray}
With abuse of language, let us mix intrinsic notation and expressions in local coordinates.
Thus we shall provisionally write $(x,\eta)$ for the cotangent vector in $T^*X$ with base $x\in X$
and fiber coordinates $\eta$ in the chosen system of Heisenberg local coordinates.
The wave front of $U^{(2)}(\tau+\tau_0)$ is contained in
$$
\left\{\left(\phi_{\tau+\tau_0}^{T^*X}\left(x,\eta\right),
\left(x,-\eta\right)\right):x\in X_1,\,\eta\in S_-\right\}.
$$
By construction, the intersection of the locus $\big\{(x,\eta):x\in X_1,\,\eta\in S_-\big\}$
with the unit sphere bundle has positive distance from
$(x_0,\alpha_{x_0})$; therefore, perhaps after restricting $X_1$, it has positive distance from $\Sigma$.
Therefore $U^{(2)}\circ \Pi$ is $\mathcal{C}^\infty$ on $(\tau_0-\epsilon,\tau_0+\epsilon)
\times X\times X$, so $\chi\cdot U^{(2)}\circ \Pi$ is $\mathcal{C}^\infty$ on $\mathbb{R}
\times X\times X$. Thus its Fourier transform is rapidly decreasing, as claimed.

Hence $S_{\chi\,e^{-i\lambda(\cdot)}}(x,x)\sim S_{\chi\,e^{-i\lambda(\cdot)}}^{(1)}(x,x)$
for $x\in X_1$ as $\lambda\rightarrow +\infty$.

Let
$F\in \mathcal{C}^\infty_0\big((0,+\infty)\big)$ be identically
$1$ in $\big(1/C_2, C_2\big)$ for some $C_2\gg 0$. Returning to
(\ref{eqn:1st-oscillatory-expression}), let us incorporate $\gamma_1\big(\eta/\|\eta\|\big)$ in the definition of
$S_{\chi\,e^{-i\lambda(\cdot)}}^{(1)}(x,x)$ in the amplitude $B$; in distributional short-hand,
\begin{eqnarray}
\label{eqn:1st-oscillatory-expression-split}
\lefteqn{S_{\chi\,e^{-i\lambda(\cdot)}}(x,x)}\\
&\sim&\frac{1}{(2\pi)^{2\mathrm{d}+1}}\,\int_{X_1}\left[\int_{\mathbb{R}^{2\mathrm{d}+1}}\int _{-\epsilon}^{\epsilon}e^{i\Psi(\tau,x,z,\eta,\lambda)}
B(\tau,x,z,\eta)\,F\left(\frac{\|\eta\|}{\lambda}\right)\,d\tau\,d\eta\right]\cdot\Pi(z,x)\,d\mu_X(z)\nonumber\\
&&+\frac{1}{(2\pi)^{2\mathrm{d}+1}}\,\int_{X_1}\left\{\int_{\mathbb{R}^{2\mathrm{d}+1}}\int _{-\epsilon}^{\epsilon}e^{i\Psi(\tau,x,z,\eta,\lambda)}
B(\tau,x,z,\eta)\,\left[1-F\left(\frac{\|\eta\|}{\lambda}\right)\right]\,d\tau\,d\eta\right\}\nonumber\\
&&\,\,\,\,\,\,\,\,\,\,\,\,\,\,\,\,\,\,\,\,\,\,\,\,\,\,\,\,\,\,\,\,\,\,\,\,\,\,\cdot\Pi(z,x)\,d\mu_X(z).\nonumber
\end{eqnarray}
Now
$|\partial_\tau\Psi|\ge C_3\big(\|\eta\|+\lambda\big)$ for some $C_3>0$ where
$1-F(\eta/\lambda)\neq 0$; therefore the second summand
is $O\left(\lambda^{-\infty}\right)$ as $\lambda\rightarrow +\infty$ uniformly in $x\in X$.

%
%
Inserting in (\ref{eqn:1st-oscillatory-expression})
the description of $\Pi$ as an FIO we get for $x\in X_1$
\begin{eqnarray}
\label{eqn:2nd-oscillatory-expression-no-heis}
\lefteqn{S_{\chi\,e^{-i\lambda(\cdot)}}\left(x,x\right)\sim}\\
&&\frac{1}{(2\pi)^{2\mathrm{d}+1}}\,\,\int_{X_1}\int_0^{+\infty}\int _{-\epsilon}^{\epsilon}\int_{\mathbb{R}^{2\mathrm{d}+1}}
e^{i\Phi_1\left(x,z,t,\tau,\eta,\lambda\right)}\cdot A\left(x,z,t,\tau,\eta,\lambda\right)\,d\mu_X(z)\,dt\,
\,d\tau\,d\eta,\nonumber
\end{eqnarray}
where
\begin{eqnarray}
\label{eqn:defn-fi-1-no-heis}
\Phi_1\left(x,\mathbf{u},z,t,\tau,\eta,\lambda\right)=\varphi\left(\tau+\tau_0,x,\eta\right)-
z\cdot \eta+t\,\psi\left(z,x\right)-\lambda\,(\tau+\tau_0),
\end{eqnarray}
and
\begin{eqnarray}
\label{eqn:defn-A}
\lefteqn{A\left(x,z,t,\tau,\eta,\lambda\right)=:}\\
&&\chi(\tau+\tau_0)\,
\varrho(z)\,a\left(\tau+\tau_0,x,z,\eta\right)\,
s\left(t,z,x\right)\,\gamma_1\left(\frac{\eta}{\|\eta\|}\right)\,
F\left(\frac{\|\eta\|}{\lambda}\right).
\nonumber
\end{eqnarray}
Apply the change of variables $t\rightsquigarrow \lambda \,t$, $\eta\rightsquigarrow \lambda\,r\,\omega$, where
$r>0$ and $\omega\in S^{2\mathrm{d}}\subseteq \mathbb{R}^{2\mathrm{d}+1}\cong \mathbb{R}\times \mathbb{R}^{2\mathrm{d}}$;
thus $\omega=(\omega_0,\omega_1)$, where $\omega_0\in \mathbb{R}$, $\omega_1\in \mathbb{R}^{2\mathrm{d}}$ and
$\omega_0^2+\|\omega_1\|^2=1$.
We can rewrite (\ref{eqn:2nd-oscillatory-expression-no-heis}) as
\begin{eqnarray}
\label{eqn:3rd-oscillatory-expression-no-heis}
\lefteqn{S_{\chi\,e^{-i\lambda(\cdot)}}\left(x,x\right)\sim
2\pi\,e^{-i\lambda\tau_0}\,\left(\frac{\lambda}{2\pi}\right)^{2\mathrm{d}+2}}\\
&&\cdot\int_{X_1}\int_0^{+\infty}\int _{-\epsilon}^{\epsilon}\int_0^{+\infty}\int_{S^{2\mathrm{d}}}
e^{i\lambda\Phi_2\left(x,z,t,\tau,r,\omega\right)}\cdot A\left(x,z,\lambda t,\tau,\lambda\,r\,\omega,\lambda\right)
\nonumber \\
&&
\cdot r^{2\mathrm{d}}\,d\mu_X(z)\,dt\,
\,d\tau\,dr\,d\omega,\nonumber
\end{eqnarray}
where
\begin{eqnarray}
\label{eqn:defn-fi-2-no-heis}
\Phi_2\left(x,z,t,\tau,r,\omega,\lambda\right)=:r\,\Big[\varphi\left(\tau+\tau_0,x,\omega\right)-
z\cdot \omega\Big]+t\,\psi\left(z,x\right)-\tau.
\end{eqnarray}

\begin{lem}
\label{lem:generating-function}
Perhaps after further replacing $S_1$ and $X_1$ with smaller open neighborhoods of $(1,\mathbf{0})\in S^{2\mathrm{d}}$ and
$x_0\in X$, respectively, we have
$$
\varphi\left(\tau+\tau_0,x,\omega\right)=\phi^X_{-(\tau+\tau_0)}(x)\cdot \omega
$$
for $x\in X_1$, $\omega\in S_1$, $|\tau|<\epsilon$.
\end{lem}

\textit{Proof.} For any $\tau$, $\varphi\left(\tau,x,\omega\right)$ is the local generating function for
the Hamiltonian flow at time $-\tau$ of the principal symbol $q$ of $Q$ on $T^*X$. In Heisenberg local coordinates,
pairs $(x,\omega)\in X_1\times S_1$ belong to a conic open neighborhood of $(x_0,\alpha_{x_0})\in \Sigma$.
On a certain conic neighborhood of $\Sigma$, on the other hand, $Q$ is microlocally equivalent to
$i\widetilde{\upsilon}_f$, hence $q$ equals the symbol of $i\widetilde{\upsilon}_f$ there.
On the same conic neighborhood of $\Sigma$, therefore, the Hamiltonian flow of $q$ is the Hamiltonian flow of
the symbol of $i\widetilde{\upsilon}_f$. However, the latter is the cotangent lift of the flow
of $\widetilde{\upsilon}_f$ on $X$.

\hfill Q.E.D.

\begin{cor}
\label{cor:nuova-phi-2}
In the same range,
\begin{eqnarray*}
\Phi_2\left(x,z,t,\tau,r,\omega,\lambda\right)=r\,\Big(\phi^X_{-(\tau+\tau_0)}(x)-
z\Big)\cdot \omega+t\,\psi\left(z,x\right)-\tau.
\end{eqnarray*}
\end{cor}

\begin{rem}
Before proceeding we note the following:

\begin{itemize}

             \item In view of the radial factor $F(r)$, integration in $dr$ is now restricted to the interval
$1/C_4\le r\le C_4$ for some $C_4\gg 0$; integration in $d\omega$ is restricted to the open neighborhood
$S_1\subseteq S^{2\mathrm{d}}$ of $(1,\mathbf{0})$.
             \item Arguing as in Lemmata 2.3 and 2.4 of \cite{p-weyl}, upon introducing a cut-off in $t$ we can further restrict integration
in $dt$ to an interval of the form $1/C_5\le t\le C_5$ for some $C_5\gg 0$.
           \end{itemize}
\end{rem}

We can now prove statement 2 in the Theorem.

For fixed $C>0$, as $\lambda\rightarrow +\infty$ we consider the family of loci
$V_\lambda\subseteq X$ given by the set of all $x\in X$ such that
$\mathrm{dist}_X\big(x,\mathrm{Fix}\left(\phi^X_{\tau_0}\right)\big)
>C\lambda^{-7/18}$.
Let $X_0\subseteq X$ be the connected component of $\mathrm{Fix}\left(\phi^X_{\tau_0}\right)$
through $x_0$, and set $F_0=:\pi(X_0)$; thus $F_0$ is the connected component of $\mathrm{Fix}\left(\phi^M_{\tau_0}\right)$
through $m_0=\pi(x_0)$, and $X_0=\pi^{-1}(F_0)$ is an $S^1$-bundle over $F_0$. Since $\pi$ is a Riemannian
submersion, we also have $\mathrm{dist}_M\big(m,\mathrm{Fix}\left(\phi^M_{\tau_0}\right)\big)
>C\lambda^{-7/18}$ if $m=\pi(x)$. As $F_0$ is a fixed clean locus for $\phi^M_{\tau_0}$,
perhaps after replacing
$C$ with a smaller positive constant we may assume that in the same range we have
$$
\mathrm{dist}_X\big(x,\phi^X_{\tau_0}(x)\big)
\ge \mathrm{dist}_M\big(m,\phi^M_{\tau_0}(m)\big)>C\lambda^{-7/18}.
$$

\begin{lem}
\label{lem:also-for-tau}
Perhaps after further decreasing $C$ and $\epsilon$,
we also have (with $m=\pi(x)$)
$$
\mathrm{dist}_X\Big(x,\phi^X_{\tau}(x)\Big)\ge \mathrm{dist}_M\big(m,\phi^M_{\tau}(m)\big)
>C\lambda^{-7/18},
$$
for all $x\in V_\lambda$ and $\tau\in (\tau_0-\epsilon,\tau_0+\epsilon)$.
\end{lem}

\textit{Proof.}
The statement of the Lemma is intrinsic, but for the proof it is convenient to make a
specific choice of Hesenberg local coordinates.

The chosen Heisenberg local coordinates centered at $x_0$ imply the prior
choice of a system of preferred local coordinates on $M$ centered at $m_0=:\pi(x_0)$,
where $\pi:X\rightarrow M$ is the projection.

Let $\mathrm{r}$ be the complex dimension of the complex submanifold $F_0\subseteq M$, and
$\mathrm{s}=:\mathrm{d}-\mathrm{r}$ its codimension.
For $m\in F_0$ let $T_m\subseteq T_mM$ and $N_m=:T_m^\perp\subseteq T_mM$ be the tangent and
normal subspaces of $F_0$ at $m$, respectively.
Furthermore, let $\exp_{m_0}^{\mathrm{tg}}:T_{m_0}\rightarrow F_0$
be the exponential map at $m_0$ of the Riemannian manifold
$F_0$, and for any $m\in F_0$ let $\exp_m^{\mathrm{nor}}:N_m\rightarrow M$ be the restriction
to $N_m$ of the exponential map of $M$.

Let us choose an orthonormal basis of $T_{m_0}$, so as to unitarily identify
$T_{m_0}\cong \mathbb{C}^{\mathrm{r}}$; also, let us choose a local unitary trivialization of the normal bundle
of $F_0$ in the neighborhood of $m_0$, so as to have smoothly varying unitarily isomorphisms $N_m\cong \mathbb{C}^{\mathrm{s}}$,
 for $m\in F_0$ near $m_0$.
In particular, by taking the normal component of the differential we have an induced unitary isomorphism
$$
d_{m_0}\phi^M_{\tau_0}:\mathbb{C}^\mathrm{s}\rightarrow \mathbb{C}^\mathrm{s},
$$
with no eigenvalue equal to $1$.

 With these isomorphisms implicit, a preferred local coordinate chart on $M$ centered at $m_0$, and defined on
some open ball centered at the origin in
$\mathbb{C}^\mathrm{r}\times \mathbb{C}^\mathrm{s}\cong \mathbb{C}^\mathrm{d}$
may then be taken
$$
\zeta(\mathbf{a},\mathbf{b})=:\exp^{\mathrm{nor}}_{\exp^{tg}_{m_0}(\mathbf{a})}(\mathbf{b}).
$$
In this chart $T_{m_0}\cong \mathbb{C}^\mathrm{r}\oplus \{0\}$ and
$N_{m_0}\cong \{0\}\oplus \mathbb{C}^\mathrm{s}$ unitarily.
Let us set $m_0+(\mathbf{a},\mathbf{b})=:\zeta(\mathbf{a},\mathbf{b})$.
Since the flow leaves $F_0$ invariant,
$$
\phi_{\tau+\tau_0}^M(m_0)=\phi_\tau^M(m_0)=m_0+\big(\upsilon(\tau),\mathbf{0}\big),
$$  where $\upsilon(\tau)=\tau \,\upsilon_f(m_0)+o(\tau)\in \mathbb{C}^\mathrm{r}$;
similarly, if $\mathbf{b}\in \mathbb{C}^\mathrm{s}$ then
clearly $\exp_{m_0}^{\mathrm{nor}}(\mathbf{b})=m_0+(\mathbf{0},\mathbf{b})$.
On the other hand, since $\phi^M_{\tau_0}$ preserves the normal geodesics at $m_0$,
we have
$$
\phi^M_{\tau_0}\big(m_0+(\mathbf{0},\mathbf{b})\big)=
m_0+\left(\mathbf{0},d_{m_0}\phi^M_{\tau_0}(\mathbf{b})\right).
$$
Therefore, for $\tau$ close to $0$ we have
$$
\phi_{\tau+\tau_0}^M(m_0+\mathbf{b})=
m_0+\Big(\upsilon(\tau),d_{m_0}\phi^M_{\tau_0}(\mathbf{b})\Big)
+
O\big(\tau\,\mathbf{b}\big).
$$

Now the previous construction may be locally smoothly deformed with $m_0\in F_0$;
furthermore, $\mathrm{Fix}\left(\phi^X_{\tau_0}\right)=\pi^{-1}(F_0)$.
If $\mathrm{dist}_X\left(x,\phi^X_{\tau}(x)\right)
>C\lambda^{-7/18}$ and $m=\pi(x)$ is in a neighborhood of $F_0$, then $m=m_0+(\mathbf{0},\mathbf{b})$ for some $m_0\in F_0$ and
$\mathbf{b}\in \mathbb{C}^\mathrm{s}$ with $\|\mathbf{b}\|>(C/2)\,\lambda^{-7/18}$.
Since a preferred local coordinate system is isometric at the origin, we get for sufficiently small
$\epsilon$ and $|\tau|<\epsilon$
\begin{eqnarray*}
\mathrm{dist}_M\left(m,\phi^M_{\tau+\tau_0}(m)\right)
&\ge& \frac 12\,\Big[\big\|\upsilon(\tau)\big\|+
\left\|\phi^M_{\tau_0}(\mathbf{b})-\mathbf{b}\right\|\Big]
+
O\big(\tau\,\mathbf{b}\big)\\
&\ge& D\,\|\mathbf{b}\|\ge D'\,\lambda^{-7/18}.
\end{eqnarray*}

Given that $\pi$ is a Riemannian submersion, and $\phi^X_{\tau'}$ covers $\phi^M_{\tau'}$, we have
$$
\mathrm{dist}_X\left(x,\phi^X_{\tau+\tau_0}(x)\right)\ge \mathrm{dist}_M\left(m,\phi^M_{\tau+\tau_0}(m)\right)
\ge D'\,\lambda^{-7/18},
$$
for some $D'>0$.

\hfill C.V.D.

\bigskip

On the way to prove statement 2, in the situation of Lemma
\ref{lem:also-for-tau} let us
rewrite (\ref{eqn:3rd-oscillatory-expression-no-heis})
as
\begin{eqnarray}
\label{eqn:4th-oscillatory-expression-no-heis}
\lefteqn{S_{\chi\,e^{-i\lambda(\cdot)}}\left(x,x\right)\sim
2\pi\,e^{-i\lambda\tau_0}\,\left(\frac{\lambda}{2\pi}\right)^{2\mathrm{d}+2}}\\
&&\cdot\int_{X_1}\left[\int_{1/D}^{D}\int _{1/D}^{D}\int_{S^{2\mathrm{d}}}\int _{-\epsilon}^{\epsilon}
e^{i\lambda\Phi_2\left(x,z,t,\tau,r,\omega,\lambda\right)}\cdot A\left(x,z,\lambda t,\tau,\lambda\,r\,\omega,\lambda\right)\right.\nonumber \\
&&\left.\cdot r^{2\mathrm{d}}\,dt\,
\,dr\,d\omega\,d\tau\right]\,d\mu_X(z),\nonumber
\end{eqnarray}
for some $D\gg 0$. Next we split the outer integral as the sum of two terms: one over
those  $z$ with
\begin{equation}
\label{eqn:1-distance-bound-xz}
\mathrm{dist}_M(z,x)>\frac C2\,\lambda^{-7/18},
\end{equation}
and one over those $z$ with
\begin{equation}
\label{eqn:2-distance-bound-xz}
\mathrm{dist}_M(z,x)\le\frac C2\,\lambda^{-7/18};
\end{equation}
here, for simplicity, we have denoted by $\mathrm{dist}_M$ the pull-back to $X\times X$
of the distance function on $M\times M$.
On the domain (\ref{eqn:1-distance-bound-xz}),
integration by parts in $t$ implies that the corresponding contribution to the asymptotics
is $O\left(\lambda^{-\infty}\right)$; in fact, a slight modification of the proof
of Lemma 2.5 of \cite{p-weyl}
shows that each integration introduces a factor $\frac 1\lambda\,\lambda^{7/9}=\lambda^{-2/9}$.
On the other hand, in view of Lemma
\ref{lem:also-for-tau} on the domain  (\ref{eqn:2-distance-bound-xz}) for $|\tau|<\epsilon$ we get
\begin{eqnarray}
\label{eqn:3-bound-distance-xz}
\lefteqn{\mathrm{dist}_X\left(\phi^X_{-(\tau+\tau_0)}(x),z\right)\ge
\mathrm{dist}_M\left(\phi^X_{-(\tau+\tau_0)}(x),z\right)}\\
&\ge& \mathrm{dist}_M\left(\phi^X_{-(\tau+\tau_0)}(x),x\right)-\mathrm{dist}_M\left(z,x\right)
 \ge \frac C2\,\lambda^{-7/18}.\nonumber
 \end{eqnarray}
Now in view of Corollary \ref{cor:nuova-phi-2}, where this holds we have
$\|\partial _\omega\Phi_2\|\ge C'\,\lambda^{-7/18}$; multidimensional integration by parts in $\omega$, therefore,
introduces at each step a factor $\frac 1\lambda\,\lambda^{7/18}=\lambda^{-11/18}$. We conclude that also in this case the contribution to
the asymptotics is $O\left(\lambda^{-\infty}\right)$.

Having established the second statement in the Theorem, let us focus on the third.
To this end, it suffices to show that there exist $a,b\in \mathbb{R}$ with $a>0$ such that for
every integer $N>0$ as $\lambda\rightarrow +\infty$ we have
\begin{eqnarray}
\label{eqn:partial-sum}
\lefteqn{S_{\chi\,e^{-i\lambda(\cdot)}}\left(x_0+\frac{\mathbf{u}}{\sqrt{\lambda}},x_0+\frac{\mathbf{u}}{\sqrt{\lambda}}\right)}\\
&\sim&
\frac{2\pi\,e^{-i\lambda\tau_0}}{f(m_0)^{\mathrm{d}+1}}\left(\frac{\lambda}{\pi}\right)^{\mathrm{d}}\,
e^{f(m_0)^{-1}\cdot\psi_2\big(d_{m_0}\phi^M_{-\tau_0}(\mathbf{u}),\mathbf{u}\big)}\,\chi(\tau_0)\cdot \left[1+\sum _{j=1}^{N}\lambda^{-j/2}G_j(x_0,\mathbf{u})\right].\nonumber\\
&&+\lambda^{\mathrm{d}}\,O\left(\lambda^{-(aN+b)}\right).\nonumber
\end{eqnarray}
Introducing Heisenberg local coordinates in (\ref{eqn:3rd-oscillatory-expression-no-heis}) we get
\begin{eqnarray}
\label{eqn:3rd-oscillatory-expression-with-heis}
\lefteqn{S_{\chi\,e^{-i\lambda(\cdot)}}\left(x_0+\frac{\mathbf{u}}{\sqrt{\lambda}},x_0+\frac{\mathbf{u}}{\sqrt{\lambda}}\right)\sim
2\pi\,e^{-i\lambda\tau_0}\,\left(\frac{\lambda}{2\pi}\right)^{2\mathrm{d}+2}}\\
&&\cdot\int_{X_1}\int_{1/D}^{D}\int_{1/D}^{D}\int _{-\epsilon}^{\epsilon}\int_{S^{2\mathrm{d}}}
e^{i\lambda\Phi_2\left(x_0+\frac{\mathbf{u}}{\sqrt{\lambda}},z,t,\tau,r,\omega
\right)}\cdot A\left(x_0+\frac{\mathbf{u}}{\sqrt{\lambda}},z,\lambda t,\tau,\lambda\,r\,\omega,\lambda\right)\nonumber \\
&&\cdot r^{2\mathrm{d}}\,d\mu_X(z)\,dt\,dr\,
\,d\tau\,d\omega.\nonumber
\end{eqnarray}

Given Lemma \ref{lem:generating-function}, for $\omega\in S_1$ and $\tau\sim 0$ we have
\begin{eqnarray}
\label{eqn:hamilton-jacobi}
\varphi\left(\tau+\tau_0,x,\omega\right)&=&\phi^X_{-(\tau+\tau_0)}
\left(x\right)\cdot \omega\\
&=&\phi^X_{-\tau_0}
\left(x\right)\cdot \omega-\tau\,\widetilde{\upsilon}_f(x)\cdot \omega+O\left(\tau ^2\right).\nonumber
\end{eqnarray}

The exponential map
$\exp_{m_0}:T_{m_0}M\rightarrow M$, composed with some unitary isomorphism $\mathbb{C}^\mathrm{d}\cong T_{m_0}M$
and restricted to some open ball $B_{2\mathrm{d}}(\mathbf{0},\delta)\subseteq \mathbb{R}^{2\mathrm{d}}\cong \mathbb{C}^\mathrm{d}$,
provides a set of preferred local coordinates for $M$ at $m_0=\pi(x_0)$.
From now on, we shall assume for convenience that these are the preferred local coordinates
underlying the chosen Heisenberg local coordinates.

Since $\phi^M_{\tau_0}(m_0)=m_0$ and $ \phi^M_{\tau_0}$ is a Riemannian isometry,
\begin{equation}
\label{eqn:hamilton-jacobi-1}
\phi^M_{-\tau_0}\left(m_0+\frac{\mathbf{u}}{\sqrt{\lambda}}\right)=m_0+\frac{1}{\sqrt{\lambda}}\,
d_{m_0}\phi^M_{-\tau_0}\left(\mathbf{u}\right).
\end{equation}
By Lemma 2.5 of \cite{p-trace} we then have
\begin{eqnarray}
\label{eqn:hamilton-jacobi-2}
\phi^X_{-\tau_0}\left(x_0+\frac{\mathbf{u}}{\sqrt{\lambda}}\right)=
\left(
\vartheta
\left(
\frac{\mathbf{u}}{\sqrt{\lambda}}
\right),
m_0+\frac{1}{\sqrt{\lambda}}\,
d_{m_0}\phi^M_{-\tau_0}
(\mathbf{u})
\right),
\end{eqnarray}
where $\vartheta$ vanishes to third order at the origin.

In Heisenberg local coordinates, $z=m_0+(\theta,\mathbf{v})$ and
$d\mu_X(z)=\mathcal{V}(\theta,\mathbf{v})\,d\theta\,d\mathbf{v}$, where $\mathcal{V}(\theta,\mathbf{0})=1/(2\pi)$.

Recalling Corollary \ref{cor:nuova-phi-2},
\begin{eqnarray}
\label{eqn:phi-2-con-diff}
\lefteqn{\Phi_2\left(x_0+\frac{\mathbf{u}}{\sqrt{\lambda}},m_0+(\theta,\mathbf{v}),t,\tau,r,\omega
\right)}\\
&=&r\left[\vartheta\left(\frac{\mathbf{u}}{\sqrt{\lambda}}
\right)-\theta\right]\,\omega_0+r\left(\frac{1}{\sqrt{\lambda}}\,
d_{m_0}\phi^M_{-\tau_0}
(\mathbf{u})-\mathbf{v}\right)\cdot \omega_1\nonumber\\
&&+t\,\psi\left(x_0+(\theta,\mathbf{v}),x_0+\frac{\mathbf{u}}{\sqrt{\lambda}}\right)-\tau.\nonumber
\end{eqnarray}

\begin{rem}
Let us pause to note the following:
\begin{itemize}
  \item If $\|\mathbf{u}\|\le C\,\lambda^{1/9}$ and $\|\mathbf{v}\|\ge 3C\,\lambda^{-7/18}$, say, then
  $\mathrm{dist}_M\big(x_0+(\theta,\mathbf{v}),x_0+\mathbf{u}/\sqrt{\lambda}\big)\ge C\,\lambda^{-7/18}$ for
  $\lambda\gg 0$; here $\mathrm{dist}_M$
  is the pull-back of the distance function on $M$.

  Given this, as in the argument in the proof of Lemma 2.5 of \cite{p-weyl} integration by parts in $dt$
  shows that the contribution to the
  asymptotics coming from the locus $\|\mathbf{v}\|\ge  3C\,\lambda^{-7/18}$ is $O\left(\lambda^{-\infty}\right)$.

  Accordingly, only a rapidly decreasing contribution is lost if an appropriate cut-off of the form
      $\gamma\left(\lambda^{7/18}\mathbf{v}\right)$ is absorbed into the amplitude,
  with $\gamma\in \mathcal{C}^\infty_0\left(\mathbb{C}^\mathrm{d}\right)$ identically equal to one near the origin.
  \item Let us perform the coordinate change $\mathbf{v}\rightsquigarrow \mathbf{v}/(r\sqrt{\lambda})$, so that integration in $d\mathbf{v}$
  is now over a ball of radius $O\left(\lambda^{1/9}\right)$. The volume element
  becomes
  $$
  d\mu_X(z)=\frac{1}{r^{2\mathrm{d}}\lambda^{\mathrm{d}}}\,\mathcal{V}\left(\theta,\frac{\mathbf{v}}{r\sqrt{\lambda}}\right)\,
  d\theta\,d\mathbf{v}.
  $$
  \item By the computations in \S 3 of \cite{sz} we have
  \begin{eqnarray*}
  \lefteqn{t\,\psi\left(x_0+\left(\theta,\frac{\mathbf{v}}{r\sqrt{\lambda}}\right),x_0+\frac{\mathbf{u}}{\sqrt{\lambda}}\right)}\\
  &=&it\left[1-e^{i\theta}\right]-\frac{it}{\lambda}\,\psi_2\left(\frac{\mathbf{v}}{r},\mathbf{u}\right)\,e^{i\theta}+
  t\,R^\psi_3\left(\frac{\mathbf{v}}{r\sqrt{\lambda}},\frac{\mathbf{u}}{\sqrt{\lambda}}\right),
  \end{eqnarray*}
  where
  $
  \psi_2(\mathbf{a},\mathbf{b})=:i\,\Im\left(\mathbf{a}\cdot \overline{\mathbf{b}}\right)-(1/2)\,\|\mathbf{a}-\mathbf{b}\|^2
  $
  and $R^\psi_3$ vanishes to third order at the origin.
\end{itemize}
\end{rem}

The upshot is that with some manipulations (\ref{eqn:3rd-oscillatory-expression-with-heis}) may be rewritten
\begin{eqnarray}
\label{eqn:4th-oscillatory-expression-with-heis}
\lefteqn{S_{\chi\,e^{-i\lambda(\cdot)}}\left(x_0+\frac{\mathbf{u}}{\sqrt{\lambda}},x_0+\frac{\mathbf{u}}{\sqrt{\lambda}}\right)\sim
\frac{e^{-i\lambda\tau_0}}{(2\pi)^{2\mathrm{d}+1}}\,\lambda^{\mathrm{d}+2}}\\
&&\cdot\int_{\mathbb{C}^\mathrm{d}}\int_{S^{2\mathrm{d}}}e^{-i\sqrt{\lambda}\mathbf{v}\cdot \omega_1}\left[\int_{-a}^a\int_{1/D}^{D}\int _{-\epsilon}^{\epsilon}\int_{1/D}^{D}
e^{i\lambda\Psi_{\mathbf{u}/\sqrt{\lambda},\omega}}\cdot B_{\mathbf{u},\mathbf{v},\omega,\lambda}\,d\theta\,dt\,d\tau\,dr\,
\right]\nonumber\\
&&\cdot\gamma\left(\lambda^{-1/6}\mathbf{v}\right)\,\,d\mathbf{v}\,d\omega,\nonumber
\end{eqnarray}
where for $\mathbf{s}\in \mathbb{C}^\mathrm{d}$, $\|\mathbf{s}\|<\delta$ we have set
\begin{eqnarray}
\label{eqn:phi-3}
\Psi_{\mathbf{s},\omega}(\theta,t,\tau,r)&=:&-r\theta\omega_0-\tau r\cdot\upsilon_f\left(x_0+\mathbf{s}\right)\cdot \omega+O\left(\tau^2\right)
\cdot r \nonumber\\
&&+it\left[1-e^{i\theta}\right]-\tau
+r\,d\Phi^X_{-\tau_0}\left(\mathbf{s}\right)\cdot \omega_1,
\end{eqnarray}
and
\begin{eqnarray}
\label{eqn:B}
\lefteqn{B_{\mathbf{u},\mathbf{v},\omega,\lambda}(\theta,t,\tau,r)=:e^{t\psi_2(\mathbf{v}/r,\mathbf{u})e^{i\theta}+it\lambda R^\psi_3\big(\mathbf{v}/(r\sqrt{\lambda}),\mathbf{u}/\sqrt{\lambda}\big)e^{i\theta}
+i\lambda r \vartheta(\mathbf{u}/\sqrt{\lambda})\omega_0}}\nonumber\\
&&
A\left(x_0+\frac{\mathbf{u}}{\sqrt{\lambda}},
x_0+\left(\theta,\frac{\mathbf{v}}{r\sqrt{\lambda}}\right),\lambda t,\tau,\lambda\,r\,\omega,\lambda\right)\, \mathcal{V}\left(\theta,\frac{\mathbf{v}}{r\sqrt{\lambda}}\right)\\
&=&e^{t\psi_2(\mathbf{v}/r,\mathbf{u})e^{i\theta}}\,\widetilde{B}_{\mathbf{u},\mathbf{v},\omega,\lambda}(\theta,t,\tau,r),
\nonumber
\end{eqnarray}
where
\begin{eqnarray*}
\lefteqn{\widetilde{B}_{\mathbf{u},\mathbf{v},\omega,\lambda}(\theta,t,\tau,r)=:
e^{i\lambda Q(\mathbf{v}/\sqrt{\lambda},\mathbf{u}/\sqrt{\lambda})}}\\
&&\cdot A\left(x_0+\frac{\mathbf{u}}{\sqrt{\lambda}},
x_0+\left(\theta,\frac{\mathbf{v}}{r\sqrt{\lambda}}\right),\lambda t,\tau,\lambda\,r\,\omega,\lambda\right)\, \mathcal{V}\left(\theta,\frac{\mathbf{v}}{r\sqrt{\lambda}}\right),
\end{eqnarray*}
with $Q$ vanishing to third order at the origin as a function of $\mathbf{v},\mathbf{u}$
(and depending on the other variables as well).
In particular, $\lambda Q(\mathbf{v}/\sqrt{\lambda},\mathbf{u}/\sqrt{\lambda})$ is bounded in our range.

We shall estimate asymptotically the inner integral by viewing it as depending parametrically on $\mathbf{v}$, $\omega_1$ and
$\mathbf{s}$, which will then be set equal to $\mathbf{u}/\sqrt{\lambda}$.
Let us
define
\begin{equation}
\label{eqn:I}
I(\mathbf{s},\mathbf{u},\mathbf{v},\omega,\lambda)=:\int_{-a}^a\int_{1/D}^{D}\int _{-\epsilon}^{\epsilon}\int_{1/D}^{D}
e^{i\lambda\Psi_{\mathbf{s},\omega}}\cdot B_{\mathbf{u},\mathbf{v},\omega,\lambda}\,d\theta\,dt\,d\tau\,dr,
\end{equation}
and
first consider the case $\mathbf{s}=\mathbf{0}$. We have
\begin{eqnarray}
\label{eqn:phi-3}
\Psi_{\mathbf{0},\omega}(\theta,t,\tau,r)=-r\theta\omega_0-\tau r\,\upsilon_f\left(x_0\right)\cdot \omega+O\left(\tau^2\right)
\cdot r +it\left[1-e^{i\theta}\right]-\tau.\nonumber
\end{eqnarray}
We have $\upsilon_f(x_0)\cdot (1,\mathbf{0})=-f(m_0)<0$, hence perhaps after restricting $S_1$
we may assume that $-\upsilon_f(x_0)\cdot \omega>\delta'>0$ for any $\omega\in S_1$.
A straightforward computation shows that

\begin{lem}
\label{lem:stationary-0}
The phase $\Psi_{\mathbf{0},\omega}$ has a unique stationary point
$$
(\theta_{\mathbf{0},\omega},t_{\mathbf{0},\omega},\tau_{\mathbf{0},\omega},r_{\mathbf{0},\omega})=
\Big(0,-\omega_0/\big(\upsilon_f(x_0)\cdot \omega\big),0,-1/\big(\upsilon_f(x_0)\cdot \omega)\Big).
$$
Furthermore, if $\Psi''_{\mathbf{0},\omega}$ is the Hessian at the critical point then
$$
\det\left(\frac{\lambda\Psi''_{\mathbf{0},\omega}
}{2\pi i}\right)=\left(\frac{\lambda}{2\pi}\right)^4\,\big(\upsilon_f(x_0)\cdot \omega\big)^2.
$$
In particular, the critical point is non-degenerate.
\end{lem}

Therefore, for $\mathbf{s}\sim \mathbf{0}$ the phase $\Psi_{\mathbf{s},\omega}$ has a unique stationary point
$
c(\mathbf{s},\omega)=(\theta_{\mathbf{s},\omega},t_{\mathbf{s},\omega},\tau_{\mathbf{s},\omega},r_{\mathbf{s},\omega}),
$
again non degenerate; furthermore, one can see by direct inspection that $c(\mathbf{s},\omega)$ is real.

By the stationary phase Lemma, for every integer $N>0$ we have
$I(\mathbf{s},\mathbf{v},\omega)=S_N(\mathbf{s},\mathbf{v},\omega)+R_N(\mathbf{s},\mathbf{v},\omega)$,
where the partial sum $S_N$ and the remainder $R_N$ are as folllows.
First,
\begin{eqnarray}
\label{eqn:expansion-general-s}
S_N(\mathbf{s},\mathbf{v},\omega)
&=&\left(\frac{2\pi}{\lambda}\right)^2\,\gamma(\mathbf{s},\omega)\,e^{i\lambda \Psi_{\mathbf{s},\omega}\big(c(\mathbf{s},\omega)\big)}
\left.\cdot\sum_{j=0}^N\lambda^{-j}\,L_j\big( B_{\mathbf{u},\mathbf{v},\omega,\lambda}\big)
\right|_{(\theta_{\mathbf{s},\omega},t_{\mathbf{s},\omega},\tau_{\mathbf{s},\omega},r_{\mathbf{s},\omega})}
\end{eqnarray}
where $\gamma(\mathbf{0},\omega)=-1/\big(\upsilon_f(x_0)\cdot \omega\big)$, $L_0$ is the identity
and any $L_j$ is a differential operator of
degree $2j$ in $\theta,t,\tau,r$, with coefficients depending on $\mathbf{s}$ and $\omega$.

On the other hand
(see Theorem 7.7.5 of \cite{h-libro}, \S 5 of \cite{sz}),
\begin{eqnarray}
\label{eqn:estimate-remainder}
\big|R_N(\mathbf{s},\mathbf{v},\omega)\big|\le \lambda^{-(\mathrm{N}+1)}C_N\,\mathrm{sup}_{|\alpha|<2N+2}\left\{\|D_{\theta,t,\tau,r}^\alpha
B_{\mathbf{u},\mathbf{v},\omega}\|\right\}.
\end{eqnarray}
In the exponent of $e^{t\psi_2(\mathbf{v}/r,\mathbf{u})e^{i\theta}+it\lambda R_3^\psi\big(\mathbf{v}/(r\sqrt{\lambda}e^{i\theta}),
\mathbf{u}/\sqrt{\lambda}\big)}$, the second summand is bounded for $\|\mathbf{u}\|,\|\mathbf{v}\|\lesssim\lambda^{1/6}$.
Therefore, (\ref{eqn:estimate-remainder}) implies

\begin{eqnarray}
\big|R_N(\mathbf{s},\mathbf{v},\omega)\big|&\le&C_N\lambda^{\mathrm{d}-3-N}\,(\|\mathbf{u}\|+\|\mathbf{v}\|)^{(4N+4)}
\,e^{-a\|\mathbf{v}-\mathbf{u}\|^2}\nonumber\\
&\le&  C_N'\lambda^{(\mathrm{d}-7/9)-(5/9)N}.              \nonumber
\end{eqnarray}
Since integration in $d\mathbf{v}$ in (\ref{eqn:4th-oscillatory-expression-with-heis}) is over a ball of radius
$O\left(\lambda^{1/9}\right)$, the overall contribution of $R_N$ is $O\left(\lambda^{2\mathrm{d}/9}\cdot \lambda^{\mathrm{d}+2}\cdot
\lambda^{(\mathrm{d}-7/9)-(5/9)N}\right)
=O\left( \lambda^{(20/9)d+11/9-(5/9)N}\right)$.

We omit the proof of the following:

\begin{lem}
\label{lem:crit-general-s}
We have $\tau_{\mathbf{s},\omega}\in \mathbb{R}$ and
$$
\tau_{\mathbf{s},\omega}=\left(d_{m_0}\phi^M_{-\tau_0}(\mathbf{s})\cdot \omega_1\right)/\big(\upsilon_f(x_0+\mathbf{s})\cdot \omega\big)+
\omega_1^tB(\mathbf{s},\omega)\omega_1,
$$
where $B(\mathbf{s},\omega)$ vanishes to second order at the origin as a function of $\mathbf{s}$.
\end{lem}

Taylor expanding $1/\big(\upsilon_f(x_0+\mathbf{s})\cdot \omega\big)$ in $\mathbf{s}$ at $\mathbf{s}=\mathbf{0}$, we get
$$
\tau_{\mathbf{s},\omega}=\left(d_{m_0}\phi^M_{-\tau_0}(\mathbf{s})\cdot \omega_1\right)/\big(\upsilon_f(x_0)\cdot \omega\big)+
\omega_1^tA(\omega,\mathbf{s})\,\omega_1,
$$
where again $A(\omega,\mathbf{s})$ vanishes to second order at the origin as a function of $\mathbf{s}$.
Consequently,
\begin{eqnarray}
\label{eqn:phase-s}
\lefteqn{i\lambda\Psi_{\mathbf{u}/\sqrt{\lambda},\omega}\left(c\left(\frac{\mathbf{u}}{\sqrt{\lambda}},\omega\right)\right)=
-i\lambda\tau_{\mathbf{u}/\sqrt{\lambda},\omega}}\\
&=&-i\sqrt{\lambda} \left(d_{m_0}\phi^M_{-\tau_0}(\mathbf{u})\cdot \omega_1\right)/\big(\upsilon_f(x_0)\cdot \omega\big)
+i \omega_1^tA_2(\mathbf{u},\omega)\omega_1+i\lambda\, P\left(\frac{\mathbf{u}}{\sqrt{\lambda}},\omega\right),\nonumber
\end{eqnarray}
where $A_2$ collects the second order terms in $A$ (as functions of $\mathbf{s}$), while
$P(\mathbf{s},\omega)$ vanishes to third order at the origin $\mathbf{s}=\mathbf{0}$,
and vanishes identically for
$\omega_1=\mathbf{0}$. In particular,
$i\lambda\,P\big(\mathbf{u}/\sqrt{\lambda},\omega\big)$ remains bounded in the given range.

Given this, (\ref{eqn:4th-oscillatory-expression-with-heis}) may be rewritten
\begin{eqnarray}
\label{eqn:5th-oscillatory-expression-with-heis}
\lefteqn{S_{\chi\,e^{-i\lambda(\cdot)}}\left(x_0+\frac{\mathbf{u}}{\sqrt{\lambda}},x_0+\frac{\mathbf{u}}{\sqrt{\lambda}}\right)=
O\left(\lambda^{(20/9)d+4/9-(5/9)N}\right)
}\\
&+&\frac{e^{-i\lambda\tau_0}}{(2\pi)^{2\mathrm{d}-1}}\,\lambda^{\mathrm{d}}\cdot
\int_{\mathbb{C}^\mathrm{d}}\int_{S^{2\mathrm{d}}}e^{i\sqrt{\lambda}\Phi_\mathbf{u}(\mathbf{v},\omega)\cdot \omega_1}\cdot K_N(\lambda,\mathbf{u},\mathbf{v},\omega)
\cdot\gamma\left(\lambda^{-1/9}\mathbf{v}\right)\,\,d\mathbf{v}\,d\omega,\nonumber
\end{eqnarray}
where
\begin{equation}
\label{eqn:phase-depends-on-u}
\Phi_\mathbf{u}(\mathbf{v},\omega_1)=:
-\left[\mathbf{v}+\frac{1}{\upsilon_f(x_0)\cdot \omega}\,d_{m_0}\phi^M_{-\tau_0}(\mathbf{u})\right]\cdot \omega_1,
\end{equation}
\begin{eqnarray}
\label{eqn:K-N}
K_N\left(\lambda,\mathbf{u},\mathbf{v},\omega\right)&=:&
e^{t_\omega\psi_2(v/r_\omega,\mathbf{u})+i \omega_1^tA_2(\mathbf{u},\omega)\omega_1
}\,\gamma\left(\frac{\mathbf{u}}{\sqrt{\lambda}},\omega\right)\\
&&\cdot\,
\left.e^{i\lambda R\big(\mathbf{v}/\sqrt{\lambda},\mathbf{u}/\sqrt{\lambda},\omega\big)
}\sum_{j=0}^N\lambda^{-j}\,L_j\left( \widetilde{B}_{\mathbf{u},\mathbf{v},\omega,\lambda}\right)
\right|_{c(\mathbf{u}/\sqrt{\lambda},\omega)}.
\nonumber
\end{eqnarray}
where we have set $t_\omega=:t(\mathbf{0},\omega)$,
$r_\omega=r(\mathbf{0},\omega)$, and $R(\cdot,\cdot,\omega)$ vanishes to third order at the origin.

Setting $\mu=\sqrt{\lambda}$, we may interpret (\ref{eqn:5th-oscillatory-expression-with-heis}) as an oscillatory integral in
$\mu$ with phase $\Phi_\mathbf{u}$.

Now $\Phi_\mathbf{u}$ has a unique critical point $(\mathbf{v}_{\mathrm{cr}},\omega_{1\,\mathrm{cr}})$, given by
$$\mathbf{v}_{\mathrm{cr}}=\frac{1}{f(m_0)}\,d_{m_0}\phi^M_{-\tau_0}(\mathbf{u}),\,\,\,\,\,\,
\omega_{1\,\mathrm{cr}}=\mathbf{0},
$$ which
is also nondegenerate: if $\Phi''_\mathbf{u}$ denotes the Hessian of $\Phi_\mathbf{u}$ at this
critical point, then
$$
\det\left(\frac{\mu \Phi''_\mathbf{u}}{2\pi i}\right)=\left(\frac{\mu}{2\pi}\right)^{4\mathrm{d}}
=\frac{\lambda^{2\mathrm{d}}}{(2\pi)^{4\mathrm{d}}}.
$$
Furthermore, since $\omega$ varies in a small neighborhood of
$(1,\mathbf{0})$, $-1/\big(\upsilon_f(x_0)\cdot \omega\big)$ is close to $1/f(m_0)$. Therefore,
upon introducing a further cut-off, we may restrict integration in $d\mathbf{v}$ to a small open neighborhood of
$f(m_0)^{-1}\,d_{m_0}\phi^M_{-\tau_0}(\mathbf{u})$, for elsewhere each integration by parts in $d\omega$
introduces a factor $\lambda^{-1/2}\,\lambda^{2/9}=\lambda^{-5/18}$, given that $A_2$ is
quadratic in $\mathbf{u}$, and $\|\mathbf{u}\|\le C\lambda^{1/9}$).

Applying the stationary phase Lemma, one sees that
the second summand on the right hand side of (\ref{eqn:5th-oscillatory-expression-with-heis})
may be rewritten
\begin{equation}
\label{eqn:remainder-sum-K}
R'_N+2\pi\,e^{-i\lambda \tau_0}\cdot \sum_{j=0}^N\lambda^{-j/2}\left.\widetilde{L}_j
(K_N)\right|_{\mathbf{v}=\mathbf{v}_{\mathrm{cr}},\omega_1=\mathbf{0}}
\end{equation}
for certain operators $\widetilde{L}_j$ of degree $2j$ in $\mathbf{v},\omega_1$,
and a remainder $R_N'$ that may be estimated by a modification of (\ref{eqn:estimate-remainder}).
More precisely, since each $\omega$ derivative brings down a quadratic factor in $\omega$, we have
\begin{eqnarray*}
\big|R_N'\big|\le C_N\lambda^{\mathrm{d}}\mu^{-(N+1)}\left(\lambda^{2/9}\right)^{2N+1}\,\lambda^{\mathrm{d}}
=C_N\,\lambda^{2\mathrm{d}-5/18-N/18}.
\end{eqnarray*}

On the other hand, in view of (\ref{eqn:B}) and (\ref{eqn:K-N}), the second summand in (\ref{eqn:remainder-sum-K})
may be rearranged as a sum of terms of the form
\begin{eqnarray}
\label{eqn:general-summand}
\lefteqn{\lambda^{-k/2}e^{f(m_0)^{-1}\psi_2\big(d_{m_0}\phi^M_{-\tau_0}(\mathbf{u}),\mathbf{u}\big)}
G_k(\mathbf{u})}\\
&&\cdot\,P_k\left(x_0,\frac{\mathbf{u}}{\sqrt{\lambda}}\right)
\cdot\,e^{if(m_0)^{-1}\lambda\big[
R_3^\psi\big(d_{m_0}\phi^M_{-\tau_0}(\mathbf{u})/\sqrt{\lambda},\mathbf{u}/\sqrt{\lambda}\big)+
\vartheta\big(\mathbf{u}/\sqrt{\lambda}\big)\big]},\nonumber
\end{eqnarray}
where $k$ is an integer and $G_k$ a polynomial.

We now use the assumption that $\mathbf{u}\in T_{m_0}M_{\tau_0}^\perp$, and
that $d_{m_0}\phi^M_{\tau_0}-I$ is invertible on $T_{m_0}M_{\tau_0}^\perp$.
This implies that $e^{f(m_0)^{-1}\psi_2\big(d_{m_0}\phi^M_{-\tau_0}(\mathbf{u}),\mathbf{u}\big)}$ is bounded by
$C e^{-a\|\mathbf{u}\|^2}$ for some $C,a>0$. Therefore, if we insert the Taylor expansion for the last two factors,
we obtain an expansion for
(\ref{eqn:general-summand}) with an $N$-th step remainder
remainder bounded by (for some polynomial $T_{kN}$)
\begin{equation}
\label{eqn:remainder-taylor}
\lambda^{-(k+N)/2}T_{kN}(\mathbf{u})\,e^{-a\|\mathbf{u}\|^2}\le
C'\lambda^{-(k+N)/2}e^{-(a/2)\|\mathbf{u}\|^2}
\end{equation}
for some $C'>0$.

This implies for (\ref{eqn:remainder-sum-K}) an
asymptotic expansion in descending powers of $\lambda^{-1/2}$, as stated in the Theorem. We now focus on
the leading term.
To this end, let us recall the asymptotic expansions for the classical symbols
in the FIO's describing $\Pi$ and $U(\tau)$:
\begin{equation}
\label{eqn:asympt-exp-amplitudes}
s(x,y,t)\sim \sum_{j\ge 0}s_j(x,y)t^{\mathrm{d}-j} \,\,\,\,\mathrm{and}\,\,\,\,a(\tau,x,z,\eta)\sim \sum_{j_{\ge 0}}
a_j(\tau,x,z,\eta),
\end{equation} where $s_0(x_0,x_0)=\pi^{-\mathrm{d}}$, and $a_j$ is homogeneous of degree $-j$ in $\eta$.
Collecting homogeneous terms of the same degree we get an asymptotic expansion
\begin{eqnarray}
\label{eqn:asympt-classical-symbols}
\lefteqn{a\left(\tau_0,x_0+\frac{\mathbf{u}}{\sqrt{\lambda}},x_0+\frac{\mathbf{v}}{\sqrt{\lambda}},\lambda r_{\mathbf{u}/\sqrt{\lambda},(1:\mathbf{0})}\omega\right)\,s\left(\lambda t_{\mathbf{0},\omega},x_0+\frac{\mathbf{v}}{\sqrt{\lambda}},
x_0+\frac{\mathbf{u}}{\sqrt{\lambda}}\right)}\nonumber\\
&&\,\,\,\,\,\,\,\,\,\,\,\,\,\,\,\sim\,\,\,\,\,
\left(\frac{\lambda}{\pi}\right)^{\mathrm{d}}\frac{1}{f(m_0)^{\mathrm{d}}}\,
a_0(\tau_0,x_0,x_0,(1:\mathbf{0}))
+\sum_{j\ge 1}\lambda^{\mathrm{d}-j/2}R_j.
\end{eqnarray}

We have $\gamma\big(\mathbf{0},(1,\mathbf{0})\big)=1/f(m_0)$ in (\ref{eqn:expansion-general-s}), and
$\mathcal{V}(\theta,\mathbf{0})=1/(2\pi)$ in (\ref{eqn:B}). Therefore,
we are left with a leading term
$$
\chi(\tau_0)\,e^{-i\lambda \tau_0}\,\left(\frac{\lambda}{\pi}\right)^{\mathrm{d}}\,\frac{1}{f(m_0)^{\mathrm{d}+1}}
e^{f(m_0)^{-1}\psi_2\big(d_{m_0}\phi^M_{-\tau_0}(\mathbf{u}),\mathbf{u}\big)}
a_0\Big(\tau_0,x_0,x_0,(1:\mathbf{0})\Big).
$$

The description of the leading coefficient is completed by the following:

\begin{lem}
\label{lem:leading-term}
$a_0\big(\tau_0,x_0,x_0,(1,\mathbf{0})\big)=2\pi$.
\end{lem}

\textit{Proof.}
The operators in our construction act on half-densities through
the trivialization of the half-density bundle offered by the volume form $d\mu_X$.
By assumption, the latter is invariant under
$\widetilde{\upsilon}_f$, that is, $\mathcal{L}_{\widetilde{\upsilon}_f}(d\mu_X)=0$,
where $\mathcal{L}$ s the Lie derivative. A straightforward computation then shows that
the subprincipal symbol of $\widetilde{\upsilon}_f$, viewed as an operator on half-densities,
vanishes identically.

Now $Q$ is microlocally equivalent to $i\widetilde{\upsilon}_f$ in a conic neighborhood
of $\Sigma$, and therefore in the same neighborhood its subprincipal symbol also vanishes.
The discussion in \S 6 of \cite{dg} then implies that in a conic neighborhood
of $(x_0,\alpha_{x_0},x_0,-\alpha_{x_0})\in \Sigma$ in the wave front of
$U(\tau_0)$ the principal symbol of $U(\tau)$ equals the natural section of
$\Omega_{1/2}\otimes L$ (the tensor product of the half-density and Maslov line bundles).

By the theory in \S 4.1 of \cite{hor-FIO-I}, $U(\tau_0)$ in local coordinates near
$x_0$ has the form $I(x,y)\,\sqrt{|dx|}\,\sqrt{|dy|}$, where
$$
I(x,y)=\int _{\mathbb{R}^{2\mathrm{d}+1}}e^{i\phi(\tau_0,x,y,\eta)}\,b(x,y,\eta)\,\big|D(\phi)\big|^{1/2}\,d\eta,
$$
where
$\phi(\tau_0,x,y,\eta)=:\varphi(\tau_0,x,\eta)-y\cdot \eta$, $D(\phi)$ is the determinant in Proposition 4.1.3
of \cite{hor-FIO-I}, and
$b(\tau_0,\cdot,\cdot)=1+b'(\cdot,\cdot)$, with $b'\in S^{-1}_{\mathrm{cl}}$. Since however we are
using the expression of $U(\tau)$ in terms of the trivialization of $\Omega_{1/2}$ given by $|d\mu_X|$, we should write
this as
$$I(x,y)=U(\tau_0)(x,y)\sqrt{\mathcal{V}_X(x)}\sqrt{\mathcal{V}_X(y)}\,\sqrt{|dx|}\,\sqrt{|dy|},
$$ so that
$$
U(\tau_0)(x,y)=\int _{\mathbb{R}^{2\mathrm{d}+1}}e^{i\phi(\tau_0,x,y,\eta)}\,b(x,y,\eta)\,\big|D(\phi)\big|^{1/2}
\,\mathcal{V}_X(x)^{-1/2}\,\,\mathcal{V}_X(y)^{-1/2}\,d\eta.
$$
Therefore,
$a_0\big(x_0,x_0,(1,\mathbf{0})\big)=2\pi\,\big|D(\phi)(\tau_0,x_0,x_0,\eta)\big|^{1/2}$.

To determine $D(\phi)$, given (\ref{eqn:hamilton-jacobi}), (\ref{eqn:hamilton-jacobi-1}) and
(\ref{eqn:hamilton-jacobi-2}) we write with $x=x_0+(\theta,\mathbf{v})$:
\begin{eqnarray}
\label{eqn:phase-tau-0}
\phi   \big(\tau_0,x_0+(\theta,\mathbf{v}),y,\eta\big)&=&
\varphi\big(\tau_0,x_0+(\theta,\mathbf{v}),\eta\big)-y\cdot \eta\\
&=&\big(\theta+\vartheta(\mathbf{u})\big)\,\eta_0+d_{m_0}\phi^M_{-\tau_0}(\mathbf{v})\cdot \eta_1-y\cdot \eta.\nonumber
\end{eqnarray}
Thus we obtain
\begin{eqnarray*}
D(\phi)(\tau_0,x_0,x_0,\eta)=\det\left(
                \begin{array}{cc}
                  \phi''_{\eta\eta} & \phi''_{\eta x} \\
                  \phi''_{y\eta} & \phi''_{yx} \\
                \end{array}
              \right)=\det\left(
                       \begin{array}{cc}
                         0 & A \\
                         -I_{2\mathrm{d}+1} & 0 \\
                       \end{array}
                     \right),
\end{eqnarray*}
where $A$ is the $(1+2\mathrm{d})\times (1+2\mathrm{d})$ matrix given by
$$
A=\left(
    \begin{array}{cc}
      1 & \mathbf{0}^t \\
      \mathbf{0} & \mathrm{Jac}_{m_0}\left(\phi^M_{-\tau_0}\right) \\
    \end{array}
  \right).
$$
The latter matrix has determinant $1$, since $\phi^M_{-\tau_0}$ is a Riemannian isometry and
and fixes $m_0$.

\hfill Q.E.D.

\bigskip

We shall now prove the last statement of the Theorem.
Let us write
 \begin{eqnarray}
\label{eqn:odd-even}
S_{\chi\,e^{-i\lambda(\cdot)}}\left(x_0+\frac{\mathbf{u}}{\sqrt{\lambda}},x_0+\frac{\mathbf{u}}{\sqrt{\lambda}}\right)=
\mathfrak{E}_\lambda(x_0,\mathbf{u})+\mathfrak{O}_\lambda(x_0,\mathbf{u}),
\end{eqnarray}
where $\mathfrak{E}$ and $\mathfrak{O}$ are even and odd functions of $\mathbf{u}$, respectively.
Thus $\mathfrak{E}$ (respectively, $\mathfrak{O}$) admits an asymptotic expansion in descending powers of $\lambda^{-1/2}$,
whose coefficients are even (respectively, odd) polynomials in $\mathbf{u}$,
and the claim is that in this expansion only integral (respectively, fractional) powers of $\lambda$ occur.

To see this, recall that the presence of fractional powers of $\lambda$ in the asymptotic expansion
of the Theorem originates from applying the stationary phase
Lemma in $\mu=\sqrt{\lambda}$ in (\ref{eqn:remainder-sum-K}) and further Taylor expanding in $\mathbf{u}/\sqrt{\lambda}$
the coefficients of the result, with remainders as in (\ref{eqn:remainder-taylor}).
Now the same expansion may be obtained as follows.

First we apply Taylor expansion in (\ref{eqn:K-N}) in
$\mathbf{u}/\sqrt{\lambda}$ and $\mathbf{v}/\sqrt{\lambda}$; if $\mathbf{v}$ is sufficiently close to
$\mathbf{v}_{\mathrm{cr}}=f(m_0)^{-1}\,d_{m_0}\phi^M_{-\tau_0}(\mathbf{u})$, we have a remainder estimate
similar to (\ref{eqn:remainder-taylor}).
The general term in this expansion will be sum of contributions of the form
$\lambda^{l-(a+b)/2}\,F_{a,b}(\mathbf{u},\mathbf{v},\omega_1)$, where $l$ is an integer and $F_{a,b}(\mathbf{u},\mathbf{v},\omega_1)$ is bihomogeneous
of bidegree $(a,b)$ in $(\mathbf{u},\mathbf{v})$.

Next we use the stationary phase Lemma in $\mu$. This can be done as above applying the
operators $\widetilde{L}_j$ in $\mathbf{v},\omega_1$, and then setting $\mathbf{v}=f(m_0)^{-1}\,d_{m_0}\phi^M_{-\tau_0}(\mathbf{u})$,
$\omega_1=\mathbf{0}$.
However, it will simplify the present discussion to proceed in the following equivalent manner.
First we make the change of variable $\mathbf{v}\rightsquigarrow
\mathbf{v}+\big(\upsilon_f(x_0)\cdot \omega\big)^{-1}\,d_{m_0}\phi^M_{-\tau_0}(\mathbf{u})$,
$\omega_1\rightsquigarrow \omega_1$, which turns the phase $\Phi_\mathbf{u}$ in (\ref{eqn:phase-depends-on-u})
into the quadratic phase $-\mathbf{v}\cdot \omega_1$ and each $F_{a,b}$
in
$$
\widetilde{F}_{a,b}(\mathbf{u},\mathbf{v},\omega_1)=:
F_{a,b}\left(\mathbf{u},\mathbf{v}-\big(\upsilon_f(x_0)\cdot \omega\big)^{-1}\,d_{m_0}\phi^M_{-\tau_0}(\mathbf{u}),\omega_1\right).
$$
The new phase has the nondegenerate critical point $\mathbf{v}=\omega_1=\mathbf{0}$.
Then we apply the stationary phase Lemma in $\mu=\sqrt{\lambda}$ in the new variables. We obtain
an asymptotic expansion given by a linear combination of terms of the following form:
\begin{eqnarray*}
\mu^{-t}\lambda^{l-(a+b)/2}\cdot \frac{\partial^2}{\partial v_{c_1}\partial \omega_{d_1}}\circ \cdots \circ
\frac{\partial^2}{\partial v_{c_t}\partial \omega_{d_t}}\left.\left(e^{t_\omega\psi_2
\left(\mathbf{v}'/r_\omega,
\mathbf{u}\right)+i\omega_1^tA_2(\omega,\mathbf{u})\omega_1}\widetilde{F}_{a,b}\right)
\right|_{\mathbf{v}=\mathbf{0},\omega_1=\mathbf{0}},
\end{eqnarray*}
where
$\mathbf{v}'=:\mathbf{v}-\big(\upsilon_f(x_0)\cdot \omega\big)^{-1}\,d_{m_0}\phi^M_{-\tau_0}(\mathbf{u})$.

Now one can see that the latter expression splits as a sum of terms of the form
$\lambda^{l-(a+b+t)/2}G(\mathbf{u})e^{f(m_0)^{-1}\psi_2\big(d_{m_0}\phi^M_{-\tau_0}(\mathbf{u}),\mathbf{u}\big)}$,
where $G$ is homogeneous of degree $2k+a+b-t$ in $\mathbf{u}$, for some integer $k$.
Thus $G$ is even if and only if $a+b+t$ is even, so that only integral powers of
$\lambda$ contribute to the asymptotic expansion of $\mathfrak{E}$.

By the same token, only fractional
(non-integral) powers of
$\lambda$ contribute to the asymptotic expansion of $\mathfrak{O}$.

This completes the proof of the Theorem.

\hfill Q.E.D.

\section{Proof of Corollary \ref{cor:global-trace-formula}.}

We want to obtain a global trace formula, that is,
an asymptotic expansion for $\int_XS_{\chi\,e^{-i\lambda(\cdot)}}(x,x)\,d\mu_X(x)$. We start off
by noticing that the integral is rapidly decreasing for $\lambda\rightarrow -\infty$ by the first
statement of the Theorem, and that for $\lambda\rightarrow +\infty$
integration may be localized near $X_0$ by the second statement. Our next step
will be to insert the local expansion in the third statement of the
Theorem within
the integral, and this requires making sense of the expression $x+\mathbf{v}$ for a variable
$x\in X_{\tau_0}$. This can be done by smoothly deforming with $x$ the construction of Heisenberg local coordinates
centered at $x$ \cite{p-trace}.
Recall that $M_{\tau_0}=:\mathrm{Fix}\left(\phi^M_{\tau_0}\right)\subseteq M$ and
$X_{\tau_0}=:\mathrm{Fix}\left(\phi^X_{\tau_0}\right)\subseteq X$, so that
$X_{\tau_0}=\pi^{-1}(M_{\tau_0})$.

Let us
fix attention on a connected component $Y=X_{\tau_0 j}$ of  $X_{\tau_0}$ at a time,
and let us set $N=\pi(Y)$. Thus $N$ is a connected component of $M_{\tau_0}$
and we shall denote its complex dimension by $\mathrm{f}$. Clearly $Y=:\pi^{-1}(N)$.

Consider a finite cover $N=\bigcup_i N_{i}$ by coordinate charts
$\beta_{i}:B_{2\mathrm{f}_0}(\delta)\rightarrow N_i$, and suppose given on each $N_i$
a local section $s_i$ of $A$ of unit norm. This induces a local chart
$\widetilde{\beta}_i:(-\pi,\pi)\times B_{2\mathrm{f}_0}(\delta)\rightarrow Y_i=:\pi^{-1}(N_i)$,
given by
$\widetilde{\beta}_i(\theta,\mathbf{r})=:e^{i\theta}\cdot s_i\big(\beta_i(\mathbf{r})\big)$.

Upon choosing the $N_i$'s sufficiently small, we may assume given for each $i$ a smooth map
$\Upsilon_i:X_{i}\times (-\pi,\pi)\times B_{2\mathrm{d}}(\mathbf{0})\rightarrow X$, such that for every
$x\in X_i$ the partial map $\Upsilon_i(x,\cdot,\cdot)$ is a Heisenberg local chart for $X$.
We may as well assume, setting $x+(\theta,\mathbf{v})=:\Upsilon(x,\theta,\mathbf{v})$, that
$x+(\theta,\mathbf{v})=\big(x+(\theta,\mathbf{0})\big)+(0,\mathbf{v})$.
Let $\mathrm{c}=:\mathrm{d}-\mathrm{f}$ be the complex codimension of $N$.
In terms of the isomorphism $\mathbb{C}^\mathrm{d}\cong \mathbb{C}^{\mathrm{f}}
\oplus \mathbb{C}^{\mathrm{c}}$, we may assume in addition that
$\Upsilon_i\big(x,\theta,(\mathbf{r},\mathbf{0})\big)\in Y$
for every $\mathbf{r}\in B_{2\mathrm{f}_0}(\delta)$.

We get a coordinate chart
\begin{eqnarray*}
\lefteqn{\widetilde{\Upsilon}_i:(-\pi,\pi)\times B_{2\mathrm{f}_0}(\delta)\times B_{2\mathrm{c}_0}(\delta)\rightarrow X,}\\
&&
(\theta,\mathbf{r},\mathbf{u})\mapsto \Upsilon_i\Big(\widetilde{\beta}_i(\theta,\mathbf{r}),0,(\mathbf{0},\mathbf{u})\Big)=
\widetilde{\beta}_i(\theta,\mathbf{r})+(\mathbf{0},\mathbf{u})
\end{eqnarray*}
with range an open neighborhood $X_i\subseteq X$ of $Y_i$.
In a natural sense, $\widetilde{\Upsilon}_i$ is \lq Heisenberg in the normal direction\rq.
Composing with $\widetilde{\beta}_i^{-1}$,
we may repackage this as a diffeomorphism
\begin{eqnarray*}
\widehat{\Upsilon}_i:Y_i\times B_{2\mathrm{c}}(\delta)\rightarrow X_i,\,\,\,\,\,\,
y\mapsto y+(\mathbf{0},\mathbf{u}).
\end{eqnarray*}
Notice that the path $y+t\,(\mathbf{0},\mathbf{u})$ meets $Y_i$ at $y$ for $t=0$ with normal velocity.
Since in the following we shall only use these normal displacements, we shall simplify notation and simply write $y+\mathbf{u}$ for
$y+(\mathbf{0},\mathbf{u})$.
Let us write $\widehat{\Upsilon}_j^*\big(d\mu_X\big)=
\mathcal{U}(y,\mathbf{u})\,d\mathbf{u}\,d\mu_{Y}(y)$, where $d\mu_{Y}$
is the natural volume form on the open subset $Y_i\subseteq Y$
(induced by the form $\omega^{\mathrm{f}}/\mathrm{f}!$ on
$N$ and the connection form); by construction, $\mathcal{U}(y,\mathbf{0})=1$ for every $y\in Y_i$.
Perhaps after composing $\widehat{\Upsilon}_j$ with a suitable change of variables of the form
$\mathbf{u}'=\mathbf{u}'(\mathbf{u},y),\,y'=y$, we may further assume that $\mathcal{U}(\mathbf{x},\mathbf{u})=1$ identically.

Finally, let $\{\gamma_{k}\}$ be a partition of unity on $N$ subordinate to the open cover $\{N_k\}$;
this may also be regarded as an $S^1$-invariant partition of unity on $Y$ subordinate to the open cover
$\{Y_k\}$. Given a tubular contraction $X'=:\bigcup _kX_k\rightarrow Y$,
the $\gamma_{k}$'s may be naturally extended to a partition of unity of $X'$. We may also arrange that
$\gamma_{k}(y+\mathbf{u})=\gamma_{k}(y)$ for all $y\in Y$ and sufficiently small
$\mathbf{u}\in \mathbb{C}^\mathrm{c}$.

We may assume that the open neighborhood $X'$ of $Y$ has positive distance from the other connected components
of $X_{\tau_0}$. By the second statement of the Theorem, therefore, $S_{\chi\,e^{-i\lambda(\cdot)}}(x,x)
=O\left(\lambda^{-\infty}\right)$ on $X'$ if $x=y+\mathbf{u}\in X'$ and $\|\mathbf{u}\|\ge
C\,\lambda^{-7/18}$.
Thus, for an appropriate radial bump function $\eta$ on $\mathbb{C}^\mathrm{c}$ identically equal to
$1$ near the origin, we have
\begin{eqnarray}
\label{eqn:partial-trace-1}
\lefteqn{\int_{X'
}S_{\chi\,e^{-i\lambda(\cdot)}}(x,x)\,d\mu_X(x)=
\sum_{k}\int_{X_k}\gamma_k(x)\cdot S_{\chi\,e^{-i\lambda(\cdot)}}(x,x)\,d\mu_X(x)}\\
&\sim&\sum_{k}\int_{Y_k\times B_{2\mathrm{c}}(\delta)}
\eta\left(\lambda^{7/18}\mathbf{u}\right)\,\gamma_k(y+\mathbf{u})\cdot S_{\chi\,e^{-i\lambda(\cdot)}}(y+\mathbf{u},
y+\mathbf{u})\,\widehat{\Upsilon}_k^*\big(d\mu_X\big)\nonumber\\
&=&\sum_{k}\int_{Y_k}\gamma_k(y)\left[\int_{B_{2\mathrm{c}}(\delta)}
\eta\left(\lambda^{7/18}\mathbf{u}\right)\,\cdot S_{\chi\,e^{-i\lambda(\cdot)}}(y+\mathbf{u},
y+\mathbf{u})\,d\mathbf{u}\right]\,d\mu_{Y}(y),\nonumber
\end{eqnarray}

Let us estimate asymptotically each of the summands in the last line of (\ref{eqn:partial-trace-1}). By
the change of variables $\mathbf{u}\rightsquigarrow \mathbf{u}/\sqrt{\lambda}$, the $k$-th
summand transforms to
\begin{eqnarray}
\label{eqn:j-th-summand-rescaled}
\lambda^{-\mathrm{c}}\,\int_{Y_k}\gamma_k\left(y\right)\cdot\left[\int_{\mathbb{C}^{\mathrm{c}}}
\eta\left(\lambda^{-1/9}\mathbf{u}\right)\, S_{\chi\,e^{-i\lambda(\cdot)}}\left(y+\frac{\mathbf{u}}{\sqrt{\lambda}},
y+\frac{\mathbf{u}}{\sqrt{\lambda}}\right)\,
d\mathbf{u}\right]\,d\mu_{Y}(y).\nonumber
\end{eqnarray}

Integration in the inner integral is now over a ball of radius $O\left(\lambda^{1/9}\right)$ centered at the origin
in $\mathbb{C}^{\mathrm{c}}$.
Since the remainder at the $N$-th step in the asymptotic expansion for $S_{\chi\,e^{-i\lambda(\cdot)}}\left(y+\mathbf{u}/\sqrt{\lambda},
y+\mathbf{u}/\sqrt{\lambda}\right)$ given by the Theorem is $O\left(\lambda^{-aN}\right)$ for
some $a>0$,
the expansion may integrated term by term. Furthermore, only the even part in the expansion gives a non-vanishing
contribution, so by item 4 in the Theorem
integration yields an asymptotic expansion in descending powers of $\lambda$.

The leading order term in the resulting asymptotic expansion is determined by computing
\begin{eqnarray}
\label{eqn:compute-leading-order}
\frac{2\pi\,e^{-i\lambda\tau_0}}{f(n)^{\mathrm{d}+1}}\left(\frac{\lambda}{\pi}\right)^{\mathrm{d}}\,
\chi(\tau_0)\lambda^{-\mathrm{c}}\,\int_{\mathbb{C}^{\mathrm{c}}}e^{f(n)^{-1}\cdot
\psi_2\big(d_{n}\phi^M_{-\tau_0}(\mathbf{u}),\mathbf{u}\big)}\,\eta\left(\lambda^{-1/9}\mathbf{u}\right)\,
d\mathbf{u},
\end{eqnarray}
where $n=\pi(y)$.

Since $\mathbf{u}$ is a normal vector, and by assumption $\mathrm{id}-d_{n}\phi^M_{\tau_0}$ is invertible
on the normal space to the fixed locus, $\Re\left(\psi_2\big(d_{n}\phi^M_{\tau_0}(\mathbf{u}),\mathbf{u}\big)\right)<-c\|\mathbf{u}\|^2$
for some $c>0$ and every $\mathbf{u}\in \mathbb{C}^{\mathrm{c}}$. Therefore, only a rapidly decreasing contribution
is lost if $\eta\left(\lambda^{-1/9}\mathbf{u}\right)$ is replaced by $1$ in (\ref{eqn:compute-leading-order})
(and similarly in the lower terms).
Performing the change of variable
$\mathbf{u}= \sqrt{f(n)}\,\mathbf{v}$, as in the derivation of (64) in \cite{p-trace} we obtain
\begin{eqnarray}
\label{eqn:compute-leading-order-1}
\lefteqn{\int_{\mathbb{C}^{\mathrm{c}}}e^{f(n)^{-1}\cdot
\psi_2\big(d_{n}\phi^M_{-\tau_0}(\mathbf{u}),\mathbf{u}\big)}\,
d\mathbf{u}}\\
&=&f(n)^{\mathrm{c}}\,\int_{\mathbb{C}^{\mathrm{c}}}e^{
\psi_2\big(d_{n}\phi^M_{-\tau_0}(\mathbf{v}),\mathbf{v}\big)}\,
d\mathbf{v}\nonumber=
f(n)^{\mathrm{c}}\cdot\frac{\pi^{\mathrm{c}}}{\det\big(\mathrm{id}_{N_{n}}
-\left.d_{n}\phi^M_{-\tau_0}\right|_{N_{n}}\big)}.
\end{eqnarray}

The leading order terms of the expansion is then the integral over $N$ of
$$
\frac{2\pi\,e^{-i\lambda\tau_0}}{f(m)^{\mathrm{f}+1}}\left(\frac{\lambda}{\pi}\right)^{\mathrm{f}}\,
\frac{\chi(\tau_0)}{\det\big(\mathrm{id}_{N_{n}}-\left.d_{n}\phi^M_{-\tau_0}\right|_{N_{n}}\big)},
$$
as claimed.

The proof is completed by repeating this argument over the set of all connected components.

\hfill Q.E.D.

\end{document}